\documentclass[]{article}

\usepackage{graphicx}
\usepackage{caption}
\usepackage{subcaption}
\usepackage{algorithmic}
\usepackage{algorithm}
\usepackage{amsmath}
\usepackage{amssymb}
\usepackage{url}
\usepackage{color}

\def\N{{\mathbb{N}}}
\def\R{{\mathbb{R}}}

\def\IMF{{\textrm{IMF}}}
\def\x{{\mathbf x}}
\def\t{{\mathbf t}}
\def\dkt{{\textrm{d}^k {\mathbf t}}}

\begin{document}

\title{Hyperspectral Chemical Plume Detection Algorithms Based On Multidimensional Iterative Filtering Decomposition}

\author{Antonio Cicone\thanks{Marie Curie fellow of the INdAM, DISIM, Universit\'a degli Studi dell'Aquila, via Vetoio 1, 67100, L'Aquila, Italy.}
\and Jingfang Liu \and Haomin Zhou\thanks{School of Mathematics, Georgia Institute of Technology, 686 Cherry St., Atlanta, GA 30332, USA.}}

\maketitle

\begin{abstract}
Chemicals released in the air can be extremely dangerous for human beings and the environment. Hyperspectral images can be used to identify chemical plumes, however the task can be extremely challenging. Assuming we know a priori that some chemical plume, with a known frequency spectrum, has been photographed using a hyperspectral sensor, we can use standard techniques like the so called matched filter or adaptive cosine estimator, plus a properly chosen threshold value, to identify the position of the chemical plume. However, due to noise and sensors fault, the accurate identification of chemical pixels is not easy even in this apparently simple situation. In this paper we present a post--processing tool that, in a completely adaptive and data driven fashion, allows to improve the performance of any classification methods in identifying the boundaries of a plume. This is done using the Multidimensional Iterative Filtering (MIF) algorithm \cite{cicone2014adaptive,cicone2015multi}, which is a non--stationary signal decomposition method like the pioneering Empirical Mode Decomposition (EMD) method \cite{huang1998empirical}. Moreover, based on the MIF technique, we propose also a pre--processing method that allows to decorrelate and mean--center a hyperspectral dataset. The Cosine Similarity measure, which often fails in practice, appears to become a successful and outperforming classifier when equipped with such pre--processing method. We show some examples of the proposed methods when applied to real life problems.

\end{abstract}

\section{Introduction}

Chemical plumes resulting from natural or anthropogenic emissions in the atmosphere can be unexpected and toxic. The detection and classification of such plumes in an efficient way would reduce the risk of harmful exposures \cite{farley2007chemical,niu2013hyperspectral}.

Recent advances in the technology provided hardware capability to easily measure the reflection energy from a large region and to generate hyperspectral images that can be used to spot and classify chemical plumes. A hyperspectral image is, in fact, a picture taken over a large number of frequencies by means of special sensors. The data produced is a hypercube whose layers are 2D images each corresponding to a different frequency. Hence each pixel of the hypercube contains a spectral signature of the corresponding physical region \cite{chang2003hyperspectral}. Given a hyperspectral image of an area where a chemical plume might be present, the plume detection problem can be divided into three subproblems. Namely: observing the presence of a plume, recognizing the chemical or chemicals contained in it, and classifying the pixels of the hyperspectral image.

Solving entirely and accurately the whole plume detection problem is not an easy task. For this reason in the end of 2012 Dimitris Manolaski and its research group at MIT Lincoln Laboratory launched the Chemical--Detection Algorithm Challenge with the idea of spurring the scientific community to solve a single subproblem. In particular they were interested in the classification stage. They gave four test datasets, each of them containing a hypercube of an area where some known chemical was released in the atmosphere, the signature of the chemical released, and the ground truth classification for each pixel. Plus they provided two blind datasets where the plume position was unknown, but it was known that a specific chemical, with a given signature, was present. The goal of the challenge was to devise a classification technique able to outperform the ones available in the literature.

Motivated by this challenge we developed, first of all, a post--processing technique, based on the Multidimensional Iterative  Filtering (MIF) method \cite{cicone2014adaptive,cicone2015multi}, which improves the classification accuracy of any given classifier. The MIF algorithm allows to decompose a nonstationary signal, associated with a linear or nonlinear system, into simpler components on which a time--frequency analysis can be performed. The development of methods for the decomposition of non--linear signals is a new and fast growing research area in Signal Processing. The interest in these kind of algorithms has been sparkled by the publication in 1998 of the breakthrough paper by Huang et al. \cite{huang1998empirical}, in which the first method of its kind, the so called Empirical Mode Decomposition (EMD) method, was proposed. The advantage of using methods like EMD or MIF versus other classical techniques available in the literature for the decomposition and time--frequency analysis of a signal, like for instance the Windowed Fourier Transform or the wavelet transform, is that these newly developed methods are designed to handle naturally non--stationary signals originated by nonlinear systems like most of the real life signals. Furthermore these techniques are completely data driven, so there is no need to make any a priori assumption on the data we want to decompose.

We note that there are many other EMD--like algorithms, such as Synchrosqueezed Wavelet Transform
\cite{daubechies2011synchrosqueezed}, the sparse time--frequency representation \cite{hou2009variant,thomas2011adaptive}, the Partial Differential Equation Transform \cite{Wang2011PDETransf,Wang2012PDETransf2}, and the ensemble EMD (EEMD) \cite{wu2009ensemble,wu2009multi}. They all aim at handling non-stationary signals associated with nonlinear systems. We select MIF because it is easy to implement, ready to be applied to 2D signals, computationally fast and most importantly it has been proved to be convergent under mild conditions, at least when applied to 1D signals.

The other main result of this paper is the development of a pre--processing algorithm, also based on MIF, which allows to decorrelate and mean--center a raw hypercube allowing to use the naive Cosine Similarity (COS) measure, also known as the Spectral Angle Mapper, as a classification technique. In fact, thanks to such a pre--processing technique the COS measure produces a classification that is always comparable and sometimes outperform the ones produced by other commonly used classifiers. The advantage of coupling the COS classifier with the proposed pre--processing method is that in this way there is no more need to estimate, unlike any other commonly used classification method, the average signature and the covariance matrix of pixels not belonging to the plume.

The rest of this paper is organized as follows: in Section \ref{sec:pb} we formalize the problem and briefly review some of the most commonly used classification methods. In Section \ref{sec:MIF} we introduce the MIF algorithm as a decomposition algorithm for one and higher dimensional signals. In Section \ref{sec:algo} we present the proposed pre-- and post--processing methods based on the MIF algorithm. Their performance are shown by means of examples in Section \ref{sec:examples}.

\section{The plume detection problem and classification algorithms}\label{sec:pb}

A hyperspectral image is a three dimensional array, whose entries are real nonnegative values, produced using hyperspectral sensors. While a standard camera sensor only captures three frequencies, the so called RGB channels, a hyperspectral sensor captures several frequencies channels at once \cite{shaw2002signal}. The output of such sensors can be represented as $F\in\R^{h\,\times\,v\,\times\,d}$ where $d$ is the number of frequency channels and $h\,\times\,v$ is the number of pixels in the image: to each pixel $p_{ij}$ it corresponds a $d$-dimensional vector $s_{ij}\in\R^d$ representing its spectral signature. Assuming we are given a hyperspectral image of an area where some chemical, potentially toxic, has been released and whose signature $s_c\in\R^d$ is given, the problem under study is the classification of pixels in $F$ as belonging or not to the chemical plume. In particular we want to assign to each pixel a score value in $[0,\, 1]$ where 0 corresponds to a chemical free pixel and 1 to a pixel containing for sure the chemical we are looking for.

Among the classification algorithms that can be used to tackle this problem the most commonly used are the Matched Filter (MF), the Normalized Matched Filter (NMF) which is also known as Adaptive Cosine or Coherence Estimator (ACE), and the plain Cosine Similarity (COS) measure also known as the Spectral Angle Mapper  \cite{manolakis2002detection,Manolakis2014Long}.

Assuming we have a single chemical released in the atmosphere whose spectral signature is $s_c$, given a pixel $p$ with signature $s$ its classification score value $y(s)$ can be computed in one of the following ways:
\begin{eqnarray}
  y(s)_{COS} &=& \frac{ \left(s^T  s_c\right)^2  }{s^T s\ s_c^T s_c}\label{COS} \\
  y(s)_{MF} &=& \frac{ \left[(s-\mu_b)^T \Sigma_b^{-1} s_c\right]^2  }{ s_c^T\Sigma_b^{-1}s_c  }\label{MF} \\
  y(s)_{ACE} &=& \frac{[(s-\mu_b)^T \Sigma_b^{-1} s_c]^2  }{ s_c^T \Sigma_b^{-1} s_c (s-\mu_b)^T \Sigma_b^{-1}(s-\mu_b)}\label{ACE}
\end{eqnarray}

Where $\mu_b\in\R^d$  is the mean signature over all the pixels in $F$ which are outside the plume, called for simplicity background pixels, and $\Sigma_b\in\R^{d\ \times\ d}$ is the covariance matrix computed considering each pixel in the background as an observation and each frequency as a variable.

We observe that if we consider the whitening and mean--centering transformation
\begin{equation}\label{eq:decorrelation}
    \widetilde{s}=\Sigma_b^{-1/2}(s-\mu_b) \quad \textrm{ and } \quad \widetilde{s}_c=\Sigma_b^{-1/2}s_c
\end{equation}
then the ACE and COS classifiers produce the very same score values
\begin{equation}\label{eq:ACEeqCOS}
    y(\widetilde{s})_{ACE} = \frac{ \left(\widetilde{s}^T  \widetilde{s}_c\right)^2  }{\widetilde{s}^T \widetilde{s}\ \widetilde{s}_c^T \widetilde{s}_c} = y(\widetilde{s})_{COS}
\end{equation}

Our goal in using all these classifiers is to properly distinguish between pixels containing the known chemical and the background pixels.

Given the scores $y(s)$ of the signatures of the pixels of a hyperspectral image $F$, for each value of a threshold $\tau\in[0,\,1]$ we can classify the pixels into four groups. If the pixel does contain the chemical and it receives a score higher than the threshold $\tau$, it is counted as a true positive, while if the assigned score is lower than $\tau$ it is counted as a false negative. The other cases are, if the pixel does not contain the chemical and it receives a score lower than $\tau$, it is counted as a true negative otherwise is a false positive.

To quantify the performance of a classifier the receiver operating characteristic (ROC) curves are commonly used \cite{fawcett2006introduction}.

A ROC curve is plotted in the Cartesian plane where on the vertical axis we have the ratio of true positive out of the total number of pixels that do contain the chemical, called True Positive rate \eqref{eq:tp_rate}, and the horizontal axis represents the ratio of false positive out of the total number of pixels that do not contain the given chemical, called False Positive rate \eqref{eq:fp_rate}. Both true positive and false positive rate are computed for different values of the threshold $\tau$.

\begin{eqnarray}
  \textrm{True Positive rate  (TPr)} &=& \frac{ \textrm{True Positive} }{ \textrm{True Chemical Pixels} }\label{eq:tp_rate} \\
  \textrm{False Positive rate (FPr)} &=& \frac{ \textrm{False Positive}  }{ \textrm{True Background Pixels} }\label{eq:fp_rate}
\end{eqnarray}

The ROC curve is a non--decreasing curve that goes from $(0,0)$ to $(1,1)$ as we vary the threshold $\tau$. The classification produced using a random guess has a corresponding ROC curve which is the straight line connecting $(0,0)$ and $(1,1)$. The larger the area under the ROC curve (AUC), the better the performance of the classifier. For more details on ROC curves we refer the reader to \cite{fawcett2006introduction}.

In the next Section we introduce the Multidimensional Iterative Filtering method which we use to devise the proposed pre-- and post--processing techniques.

\section{Multidimensional Iterative Filtering (MIF)}\label{sec:MIF}

Finding features, structures or quasiperiodicities in a non--stationary signal generated by a non--linear system can be quite challenging. However many real life problems require to handle and analyze such signals. We can think, for instance, to finance or climatic studies where the identification of recurrent patterns in data collected over time would be extremely valuable. For this reason in the last two decades many tools have been devised to deal with such signals. The goal is to decompose a signal into a finite number of simpler components called intrinsic mode functions (IMFs) plus a trend. Intuitively speaking, an IMF is a function oscillating around zero. Huang et al. in \cite{huang1998empirical} defined an IMF as a function such that: in its entire domain the number of extrema and the number of zero crossings must either equal or differ at most by one; at any point, the mean value of the upper envelope defined by connecting its local maxima with splines and the lower envelope obtained by connecting its local minima is zero. From these intuitive definition it follows that sinusoidal functions are proper IMFs, but IMFs are not limited to those functions. We can have changes in the amplitude and/or the frequency as well as discontinuities in an IMF or in its derivatives.

There are many algorithm developed to decompose a signal into IMFs. The commonly used ones work either by iteration or by optimization. The pioneer work, called the Empirical Mode Decomposition (EMD) method published by Huang et al. in 1998 \cite{huang1998empirical}, has an iterative structure,
the same as the Iterative Filtering method developed by Zhou et al. \cite{lin2009iterative}. While examples of methods which are based on optimization are the Sparse Time--Frequency Representation by Hou et al. in \cite{hou2009variant,thomas2011adaptive} as well as the Synchrosqueezed Wavelet Transforms proposed by Daubechies et al. in \cite{daubechies2011synchrosqueezed}. For an overview on such techniques we refer the reader to \cite{cicone2014adaptive}.

In this paper we focus on the Multidimensional Iterative Filtering algorithm \cite{lin2009iterative,cicone2015multi} (MIF), which is an extension to higher dimensions of the Iterative Filtering method that is known to be convergent under mild assumptions, at least for 1D signals, and numerical experiments showed its stability \cite{cicone2014adaptive}. We point out here that the other known iterative method, the EMD algorithm suffers of known instabilities and its convergence is still an open problem. While a variant of the EMD method, called Ensemble EMD (EEMD), allows to solve the instability issues and it has been extended to higher dimensions \cite{wu2009multi}, its convergence remains an open problem since EEMD is based on the repeated usage of the EMD algorithm.

In the MIF algorithm\footnote{The MIF code can be found at \url{www.cicone.com}} each IMF is produced convolving iteratively the signal with a low pass filter $w(\t)$, $\t\in \Omega\subset\R^k$ like, for example, a Fokker--Planck filter \cite{cicone2014adaptive} which has the nice property of being compactly supported and smooth on its entire domain. In particular, given the signal $s(\x), \x\in\R^k$, the MIF algorithm works as follows

\begin{algorithm}
\caption*{\textbf{MIF Algorithm} IMF = MIF$(s)$}
\begin{algorithmic}
\STATE IMFs = $\left\{\right\}$
\WHILE{the average number of extrema of $s \geq 2$}
      \STATE  compute the filter support $\Omega$ for $s$
      \STATE $s_1 = s$
      \WHILE{the stopping criterion is not satisfied}
                  \STATE  $s_{n+1}(\x) = s_{n}(\x) -\int_{\Omega} s_n(\x+\t)w(\t)\dkt$
                  \STATE  $n = n+1$
      \ENDWHILE
      \STATE IMFs = IMFs$\,\cup\,  \{ s_{n}\}$
      \STATE $s=s-s_{n}$
\ENDWHILE
\STATE IMFs = IMFs$\,\cup\,  \{ s\}$
\end{algorithmic}
\end{algorithm}

Regarding the selection of the filter support $\Omega$ we want to base it adaptively on the signal. Following what proposed in \cite{lin2009iterative} for the 1D case where the length of the support is computed as $4N/K$, where $N$ is the number of sample points in the 1D signal and $K$ is the number of its extrema, we compute for each dimension of the signal the average support length. Then we use this information to either find the radius of a spherical support or to identify the radii of an ellipsoidal $\Omega\subset\R^k$.

For the stopping criterion there are several possibilities. One way of doing it, as explained for the 1D case in \cite{cicone2014adaptive,huang1998empirical}, is to consider the relative change
\begin{equation}\label{eq:SD}
SD = \frac{  \|s_{n+1}(\x) - s_{n}(\x)\|_{l^2}   }{ \| s_{n}(\x) \|_{l^2}   }.
\end{equation}
and discontinue  the inner loop as soon as this relative change is less than a given threshold. The value of SD is usually set around $0.001$ in our experiments.

The algorithm stops as soon as the residual $r(t)$ has a dimension which contains at most one extrema. So in the end the given signal is decomposed as $s(\x) = \sum_{j=1}^m \IMF_j(\x) + r(\x)$, where $m$ is the number of IMFs produced by the algorithm.

\section{Pre-- and post--processing algorithms based on MIF}\label{sec:algo}

In this section, we propose both a post-- and a pre--processing technique based on the Multidimensional Iterative Filtering (MIF) algorithm.

The post--processing method aims to increase the performance of any given classifier. Assuming we are dealing with a chemical plume which is not punctiform, so that at least a few pixels can be classified as containing a certain chemical, the performance of a classifier are effected by noise or single sensor misfunctioning \cite{Manolakis2014Long}. The idea, as pointed out by Wager and Walther in \cite{Wager2012DTRA}, is that pixels classified as containing the chemical we want to detect have, by assumption, to be surrounded by other pixels containing that very same chemical. We can improve the classification through the reinforcement of the spatial correlation. To do so, using the MIF method, we decompose the non--stationary and bidimensional signal, produced by a classifier, into high and low frequency components in an adaptive and data driven fashion. Then, we remove the high frequencies components, which do not contain information related with the plume, preserving instead the low frequency components. In doing so we reduce the rate of misclassified pixels improving the performance of any possible classification algorithm, not only the ones reviewed in this work. The choice of the MIF method is very natural since, unlike any other technique available in the literature, this method is designed to handle naturally non--stationary signals, it does not require to make any a priori assumption on the signal itself and its convergence has been proved when applied to 1D signals.

\begin{description}
\item[Post--processing method (\textbf{PostP})] Given the classification matrix $C\in\R^{h\,\times\,v}$, produced by a given classifier, we decompose $C$ using the MIF method into its first IMF $I_{1}$ and a remainder/trend $R$. We then remove the first IMF from $C$ producing the new classification $\widehat{C}=C-I_1\equiv R$. We represent this method using the operator $\textrm{PostP}:\,\R^{h\,\times\,v}\mapsto\R^{h\,\times\,v},\; C\rightarrow \widehat{C}=\textrm{PostP}(C)$.
\end{description}

We observe that, since the convergence of this method applied to the decomposition of 2D signals is still an open problem, we checked the meaningfulness of the outcome of such a post--processing using the 1D convergent MIF method . We consider, in fact, the classification scores first along each row and then along each column of the classification matrix, and vice versa. The results of these two procedures prove to be close to the one obtained post--processing directly the 2D classification using the 2D MIF method.

The other method we present here is a pre--processing technique. The goal is to devise a completely data driven decorrelation and mean--centering technique that can be applied to a hyperspectral dataset before classifying it with the plain Cosine Similarity (COS) measure. The idea behind this pre--processing method comes from the observations, we made previously, regarding the ACE and COS classifiers. In particular we know that ACE requires the a priori knowledge of the mean signature and covariance matrix of the background and we observed that if we can somehow whiten and mean--center the hyperspectral image pixel signatures then the COS and ACE classifiers produce the very same classification scores \eqref{eq:ACEeqCOS}. To achive this goal we propose the following procedure. First of all we subtract from each pixel signature $s$ the average signature $\mu$ of the entire hyperspectral dataset. This is done in order to remove possible artifacts introduced uniformly in all the pixel signatures by the sensor or the device used to capture the hyperspectral image. Then we decompose each pixel signature, using the MIF method, into IMFs and a trend. Hence we subtract such trend from the original signature. The result is a hyperspectral dataset which is mean--centered and decorrelated in a local way and it is ready to be classified using the plain Cosine Similarity (COS) measure.

\begin{description}
 \item[Pre--processing method (\textbf{PreP})] Given a hypercube $F$ and the signature $s_{ij}\in\R^d$ of the $(i,\,j)$ pixel, we first compute $\mu\in\R^d$, mean signature of all the pixels in $F$, and we subtract it from each pixel signature to produce $\widetilde{s}_{ij}=s_{ij}-\mu$. Then using the MIF method we decompose $\widetilde{s}_{ij}$ into its IMFs $\left\{I_{n}\right\}_{n=1,\ldots N}$, for some $N\in\N$, and a trend $r$ so that $\widetilde{s}_{ij} = \sum_{n=1}^N I_{n}^{(i,j)} + r$. Finally we remove from the signature $\widetilde{s}_{ij}$ the remainder/trend $r$. So each pixel signature becomes $\widehat{s}_{ij}=\sum_{n=1}^{N} I_{n}^{(i,j)}$. We can summarize the entire procedure using the operator $\textrm{PreP}:\,\R^{h\,\times\,v\,\times\,d}\mapsto\R^{h\,\times\,v\,\times\,d},\; F\rightarrow \widehat{F}=\textrm{PreP}(F)=\left[\widehat{s}_{ij}\right]_{i,\,j}$.
\end{description}

The proposed pre--processing technique can be applied to hyperspectral images. It allows, in particular, to decorrelate and mean--center the hypercubes. Since we simply decorrelate and we do not whiten the dataset we do not expect COS and ACE classifications to match in general. However from the experiments we run we observe that COS equipped with such a pre--processing method proves to have performance that are always similar and sometimes even better than the one of the ACE classifier. Therefore, from a COS classification point of view, the global and local trend removal performed in the proposed pre--processing method proves to be an equivalent procedure to the mean--centering and whitening of the data, as shown in the next section.

\section{Examples}\label{sec:examples}

In this section we show the performance of the proposed pre-- and post--processing methods when applied to the plume detection problem. We consider the datasets provided for the 2012 DTRA Chemical Detection Challenge\footnote{The datasets can be obtained from the National Science
Foundation (NSF). However, they are provided only to researchers and teams supported by NSF itself.}. They contain both real world and synthetic hyperspectral images, and the signatures of the chemical contained in the plume are always included. In this work we focus on the real world datasets, since the analysis of the synthetic ones appear to be trivial in general.
We test our methods first on a dataset where the ground truth is known, since the actual position of the chemical plume is given, and then on a blind one that have been provided for the contest. With the former we can evaluate the performance of each classifier using ROC curves.

We point out here that the ground truths provided in the datasets for the 2012 DTRA Chemical Detection Challenge assume a classification of the hyperspectral image pixels into three groups:  ``inside the plume'', ``outside the plume'' and ``close to the boundaries of the plume'', as shown, for instance, in Figure \ref{fig:dugway_rel_r134_Plume}.  As explained in Section \ref{sec:pb}, in this field of research the goal is to devise a method able to classify pixels using just two classes: pixels containing the known chemical and background pixels. The reason why the ground truth includes the additional class of pixels close to the boundaries of the plume is because classifying such pixels as either inside or outside the plume is a hard problem. Therefore, following what suggested by Manolakis and his research group for the 2012 DTRA Chemical Detection Challenge, we use the ground truth to check the performance of a classifier only on pixels classified as either ``inside the plume'' or ``outside the plume''. This means that the ROC curves presented in this Section are built using solely the pixels that are not contained in the set ``close to the boundaries of the plume''.
It is clear that if we use only ROC curves built in the aforementioned way we can compare the performance of different classifiers on pixels outside and inside the plume, but we are unable to compare them on pixels close to the boundaries of the plume. For this reason we show both ROC curves and the pixel classification values produced by each classifier in order to give readers a more comprehensive picture of the performance of each classifier.

In the blind case, instead, it is impossible to plot ROC curves since the ground truth is unknown. In this second case we show the performance of a classifier simply plotting the pixel classification values, also known as detection maps.

We plot the pixel classification values using ``flipud'' colormap, option ``hot'', in Matlab.  In all the examples presented in the following the colormap is calibrated identically.

We would like to note that we design MIF and the pre-- and post--processing techniques for hyperspectral data classifications not to compete with the existing classification methods,  but to use them as complementary strategies to further improve their performances. In this paper, we present the performances of the proposed pre-- and post--processing techniques when applied to standard classifiers like ACE, MF and COS. We point out that they can be easily adopted to other classifiers as well.

\subsection{Case of known plume position: hypercube ``Location 1 released r134a'' }

\begin{figure}[H]
        \begin{subfigure}[b]{0.48\textwidth}
                \centering
                \includegraphics[width=\textwidth]{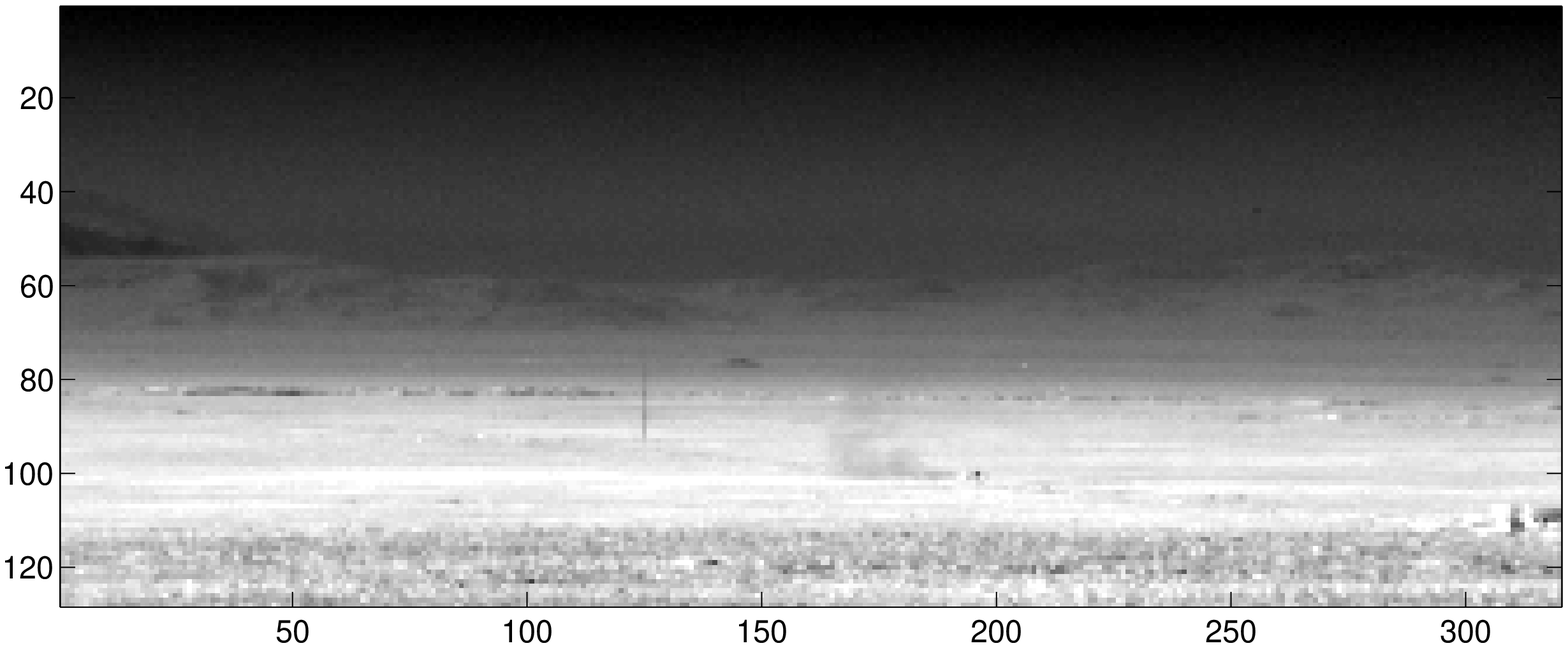}
                \caption{Contrast--enhanced spectral--mean image using the raw data from the dataset \emph{Location 1 released r134a}}
                \label{fig:dugway_rel_r134_Image}
                \end{subfigure}\hskip 2mm
                \begin{subfigure}[b]{0.48\textwidth}
                \centering
        \includegraphics[width=\textwidth]{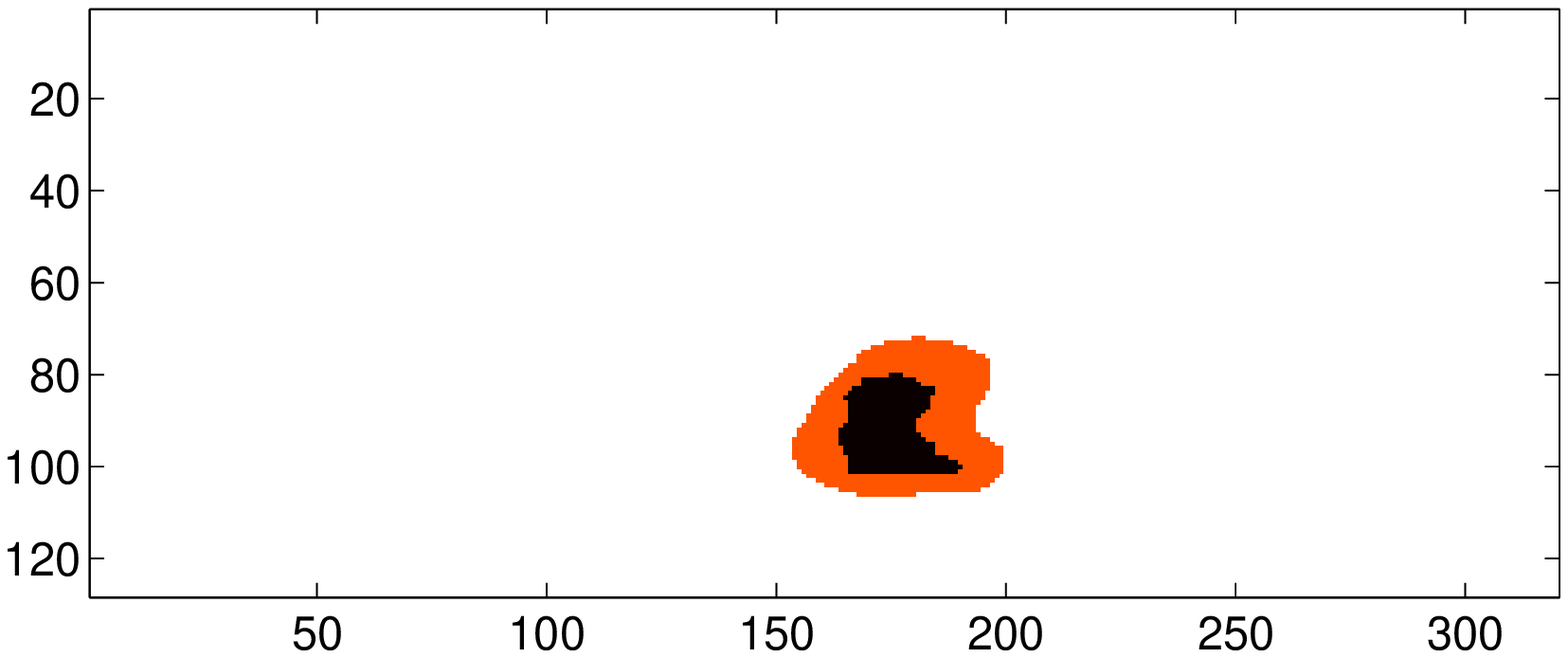}
        \caption{Ground truth classification. In black plume pixels, in orange pixels close to the boundaries of the plume.}
        \label{fig:dugway_rel_r134_Plume}
        \end{subfigure}
\caption{}\end{figure}

\begin{figure}[H]
\begin{subfigure}[b]{0.48\textwidth}
\centering
        \includegraphics[width=\textwidth]{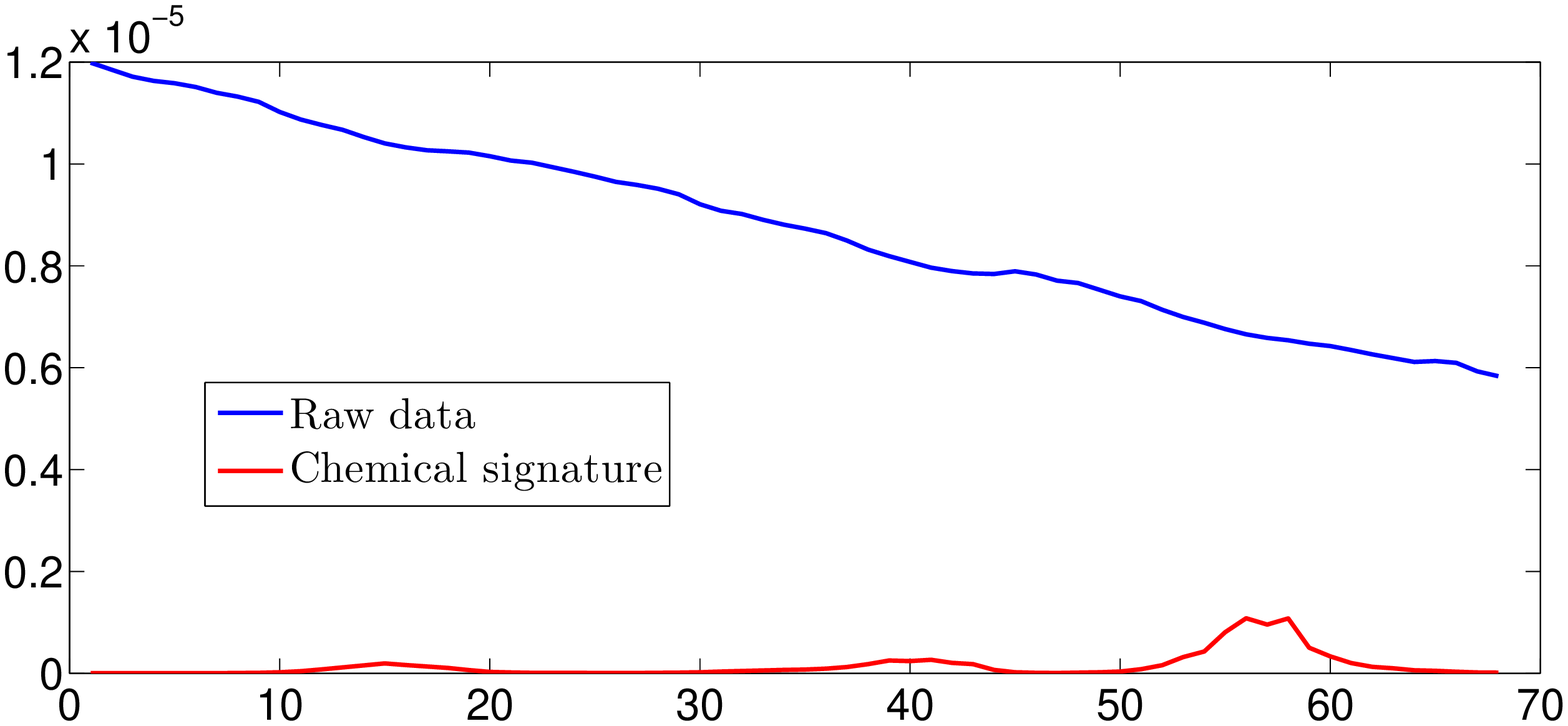}
                \caption{Raw spectral signature of a pixel inside the plume and the chemical r134 signature}
                \label{fig:dugway_rel_r134_RawPlumePixel1}
        \end{subfigure}\hskip 2mm
      \begin{subfigure}[b]{0.48\textwidth}
     \centering
     \includegraphics[width=\textwidth]{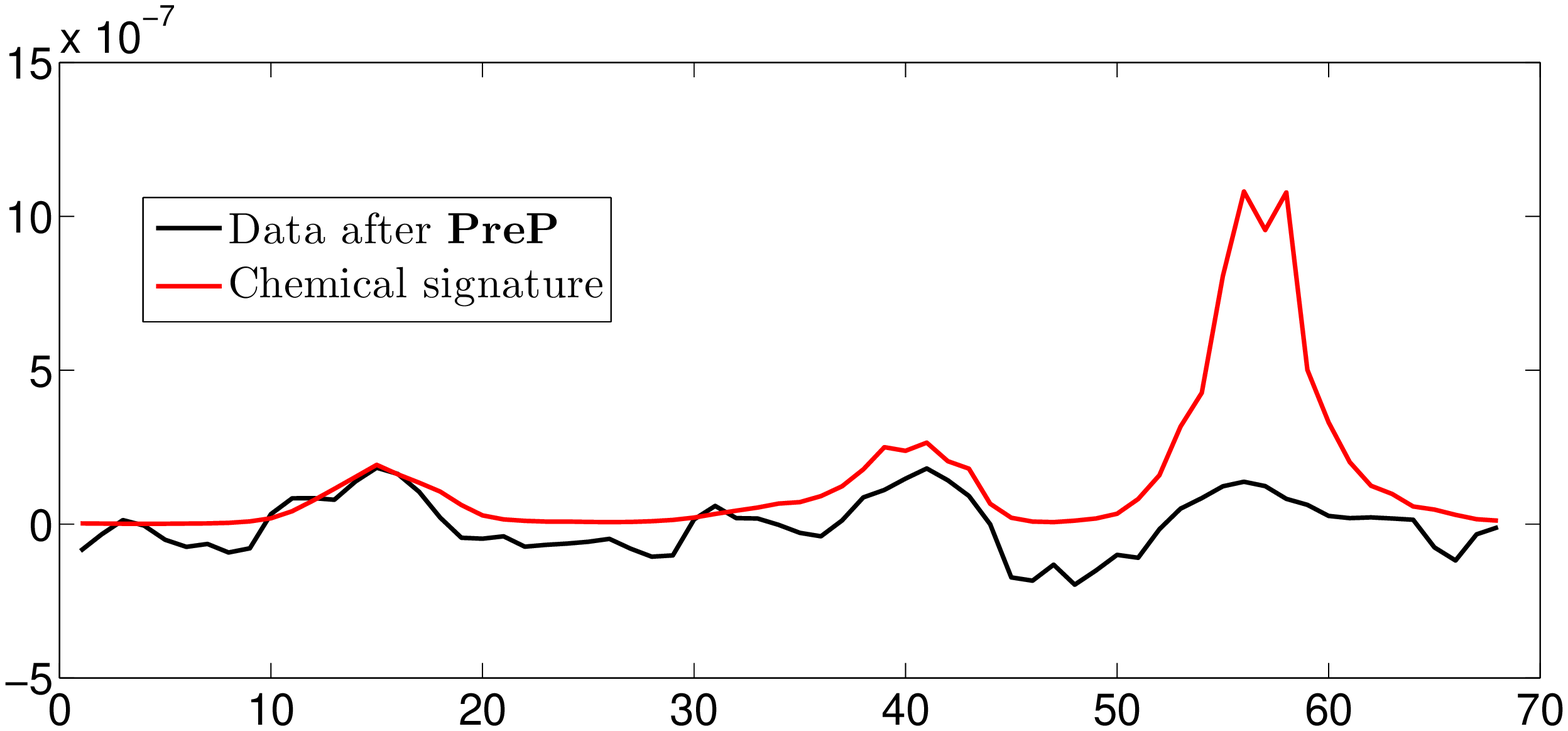}
     \caption{Pre--processed spectral signature of a pixel inside the plume}
     \label{fig:dugway_rel_r134_PreprocPlumePixel1}
     \end{subfigure}
\caption{}\end{figure}

The first dataset we test contains the hypercube plotted in Figure \ref{fig:dugway_rel_r134_Image} taken in an area where the chemical r134 has been released. The signature of this chemical is provided in the dataset and shown in red in Figure \ref{fig:dugway_rel_r134_PreprocPlumePixel1}. The given ground truth is shown in Figure \ref{fig:dugway_rel_r134_Plume} where pixels inside the plume are colored black, while ones close to the boundaries of the plume are colored orange. We observe from this figure that the boundary region, that we are excluding in the evaluation of the ROC curves, as explained previously, is pretty wide. It becomes evident that using only ROC curves and area under the curve (AUC) values provides only partial information on the true performance of each classifier.

\begin{figure}[H]
        \begin{subfigure}[b]{0.48\textwidth}
        \centering
        \includegraphics[width=\textwidth]{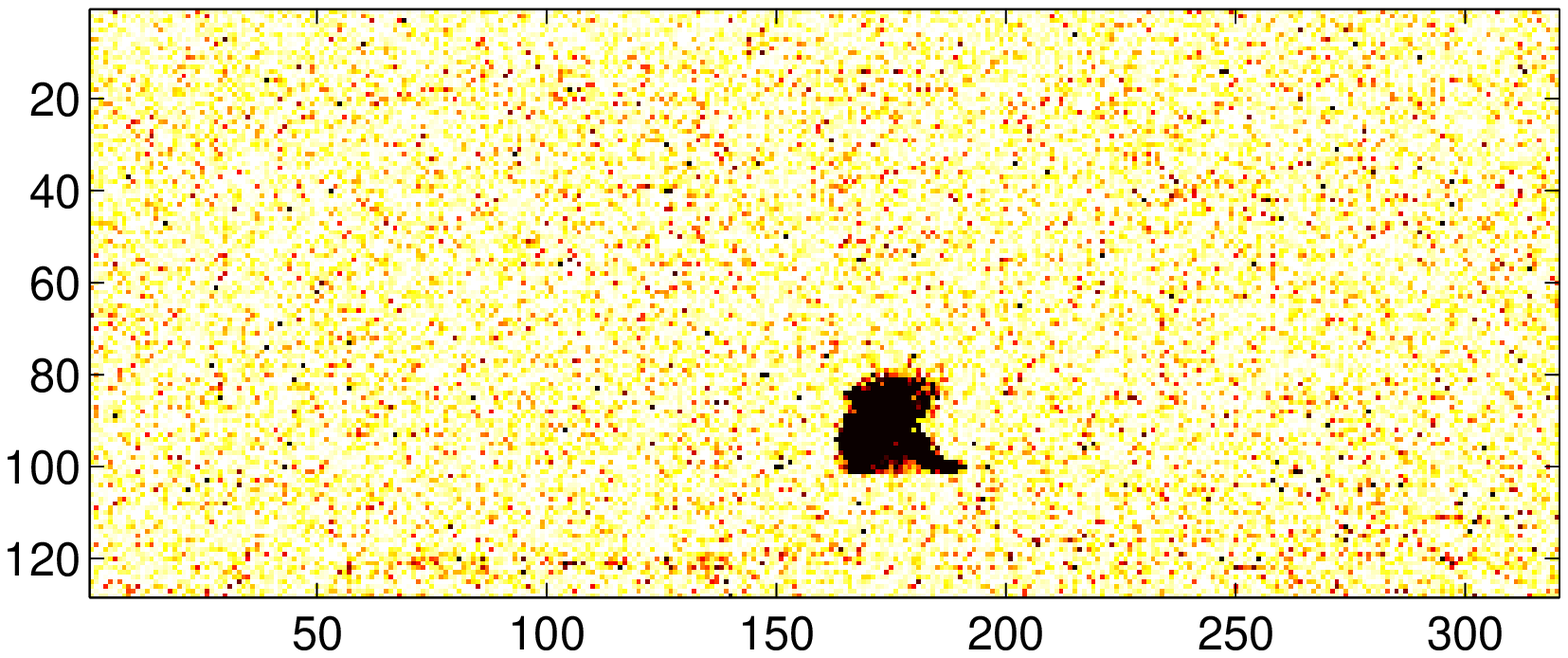}
        \caption{Detection map of the ACE classification of the raw data}
        \label{fig:dugway_rel_r134_det_ACE}
        \end{subfigure}\hskip 2mm
        \begin{subfigure}[b]{0.48\textwidth}
        \centering
        \includegraphics[width=\textwidth]{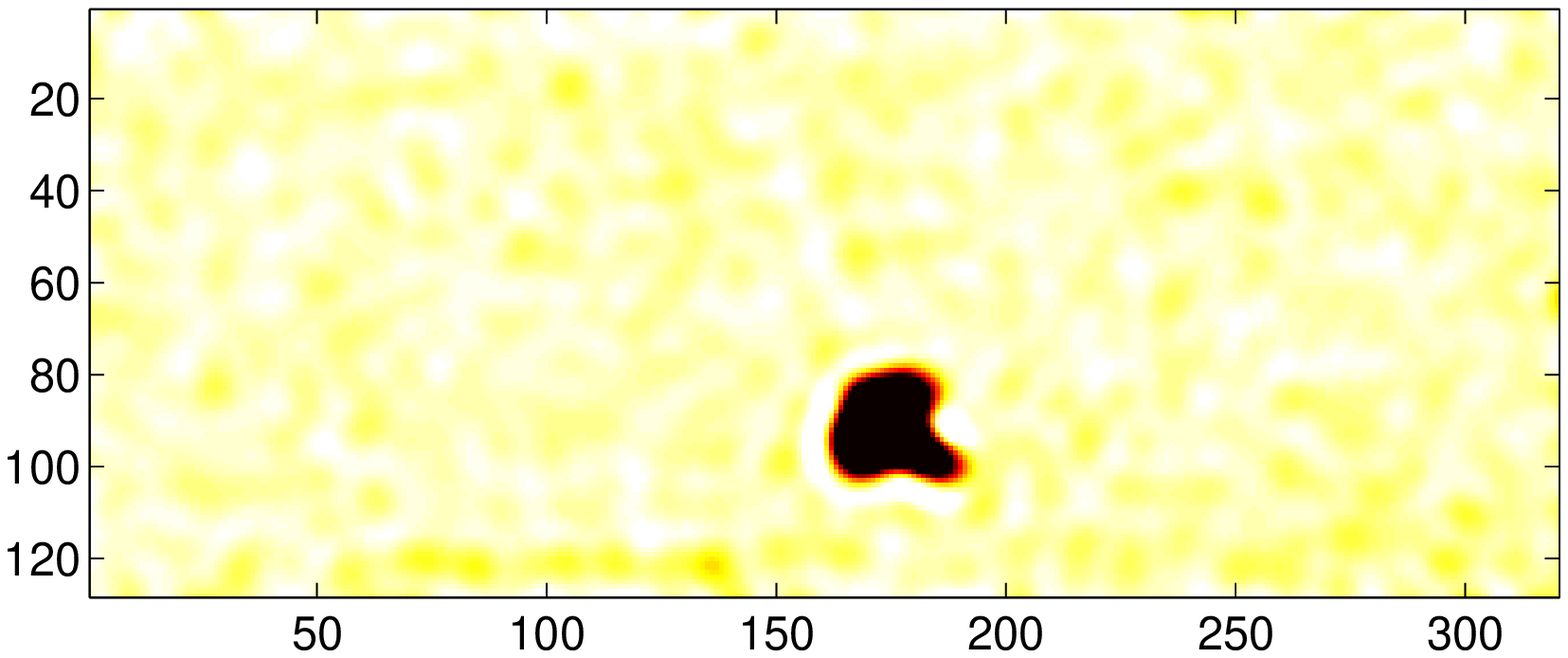}
                \caption{Detection map after applying the \textbf{PostP} method to the ACE classification of the raw data}
                \label{fig:dugway_rel_r134_det_ACE_PP}
        \end{subfigure}
\caption{}
\end{figure}

\begin{figure}[H]
\begin{subfigure}[b]{0.48\textwidth}
                \centering
                \includegraphics[width=\textwidth]{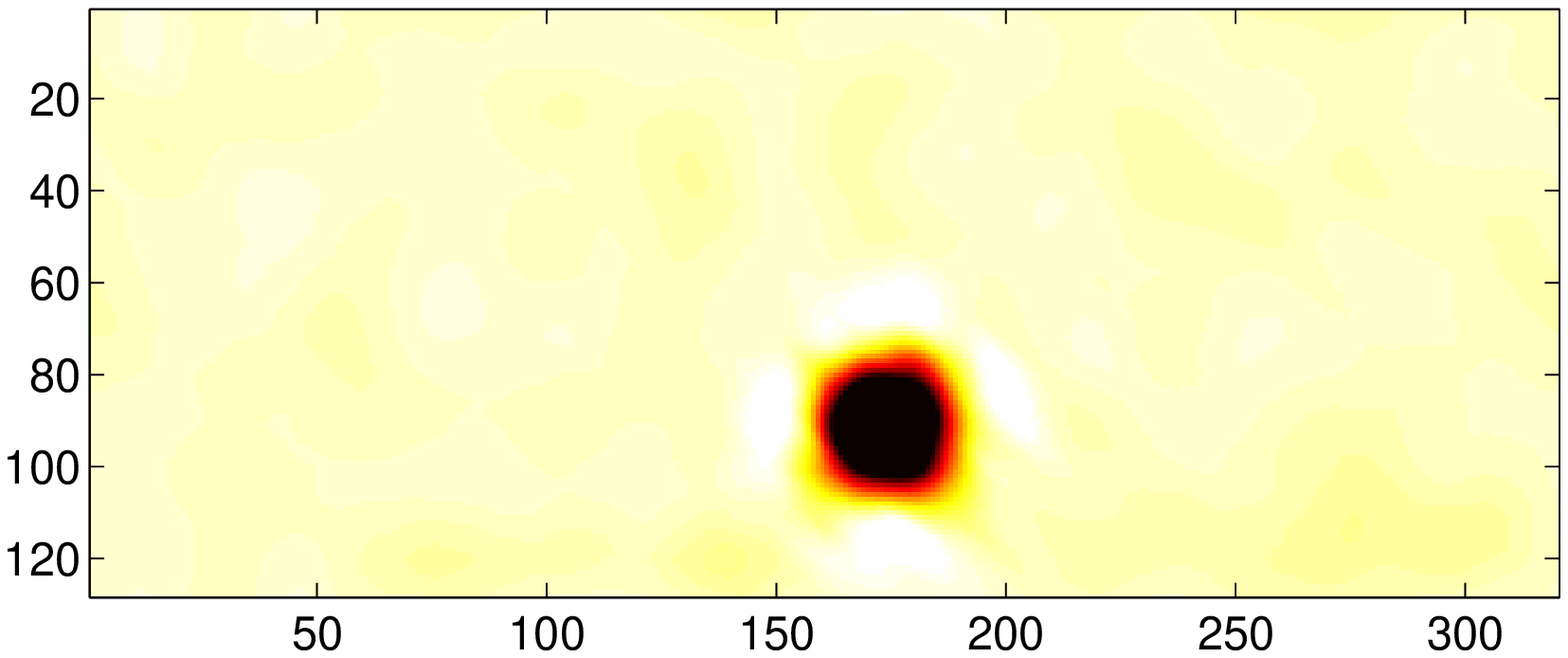}
        \caption{Detection map of the Wavelets post--processed ACE classification of the raw data}
        \label{fig:dugway_rel_r134_det_ACE_Wav_PP}
        \end{subfigure}
      \hskip 2mm
     \begin{subfigure}[b]{0.48\textwidth}
     \centering
        \includegraphics[width=\textwidth]{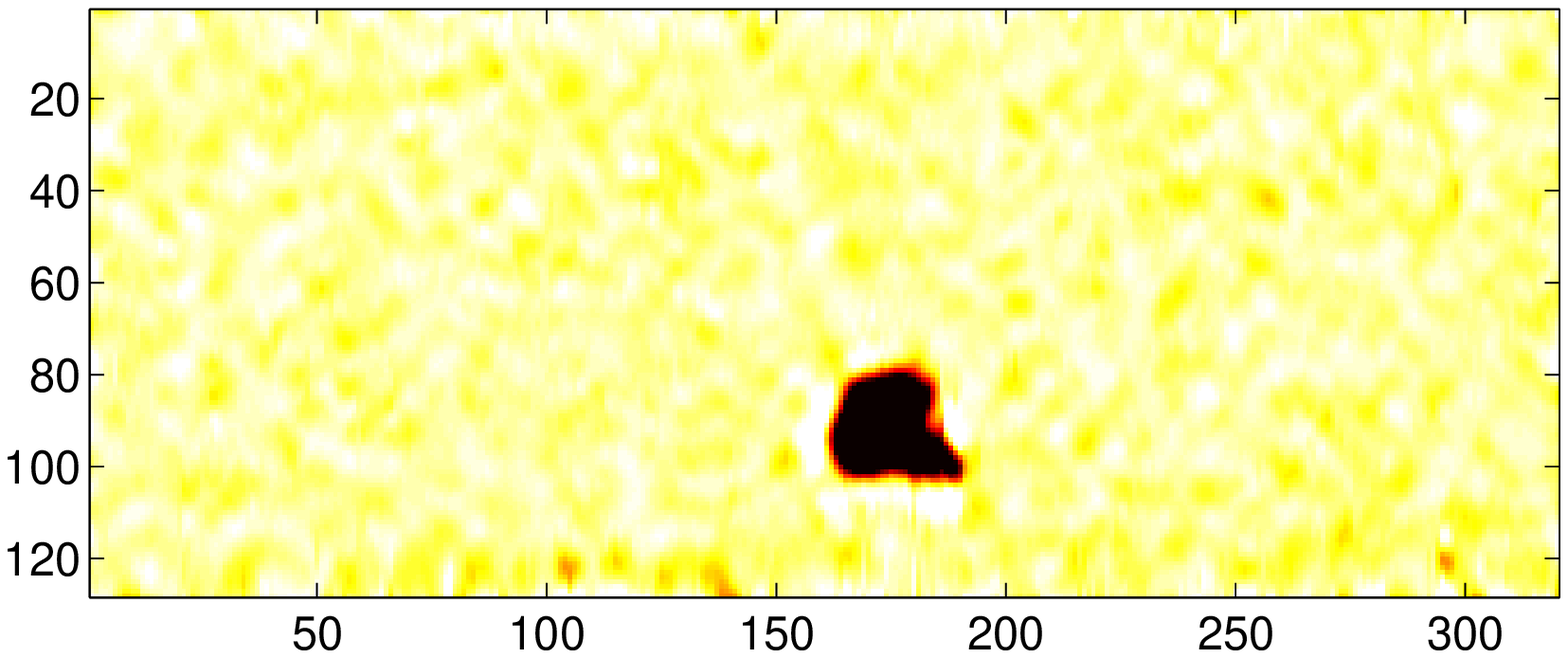}
        \caption{Detection map of the MEEMD post--processed ACE classification of the raw data}
        \label{fig:dugway_rel_r134_det_ACE_PP_MEEMD}
        \end{subfigure}
\caption{}\end{figure}

\begin{figure}[H]
        \begin{subfigure}[b]{0.48\textwidth}
     \centering
     \includegraphics[width=\textwidth]{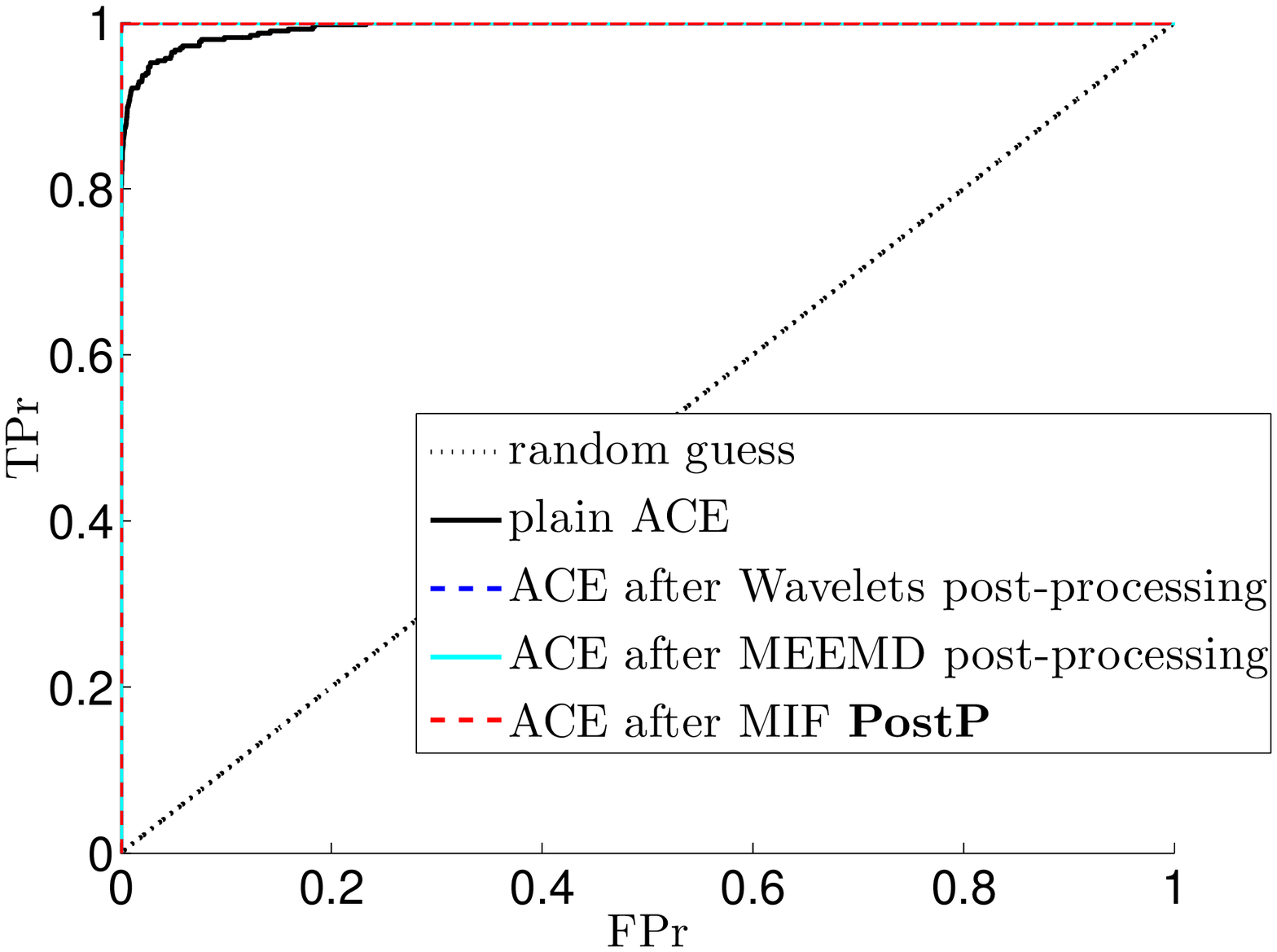}
     \caption{ROC curves  for the ACE classification of the hypercube \emph{Location 1 released r134a} with and without post--processing}
     \label{fig:dugway_rel_r134_roc_ACE_PP_noLog}
     \end{subfigure}\hskip 2mm
     \begin{subfigure}[b]{0.48\textwidth}
     \centering
     \includegraphics[width=\textwidth]{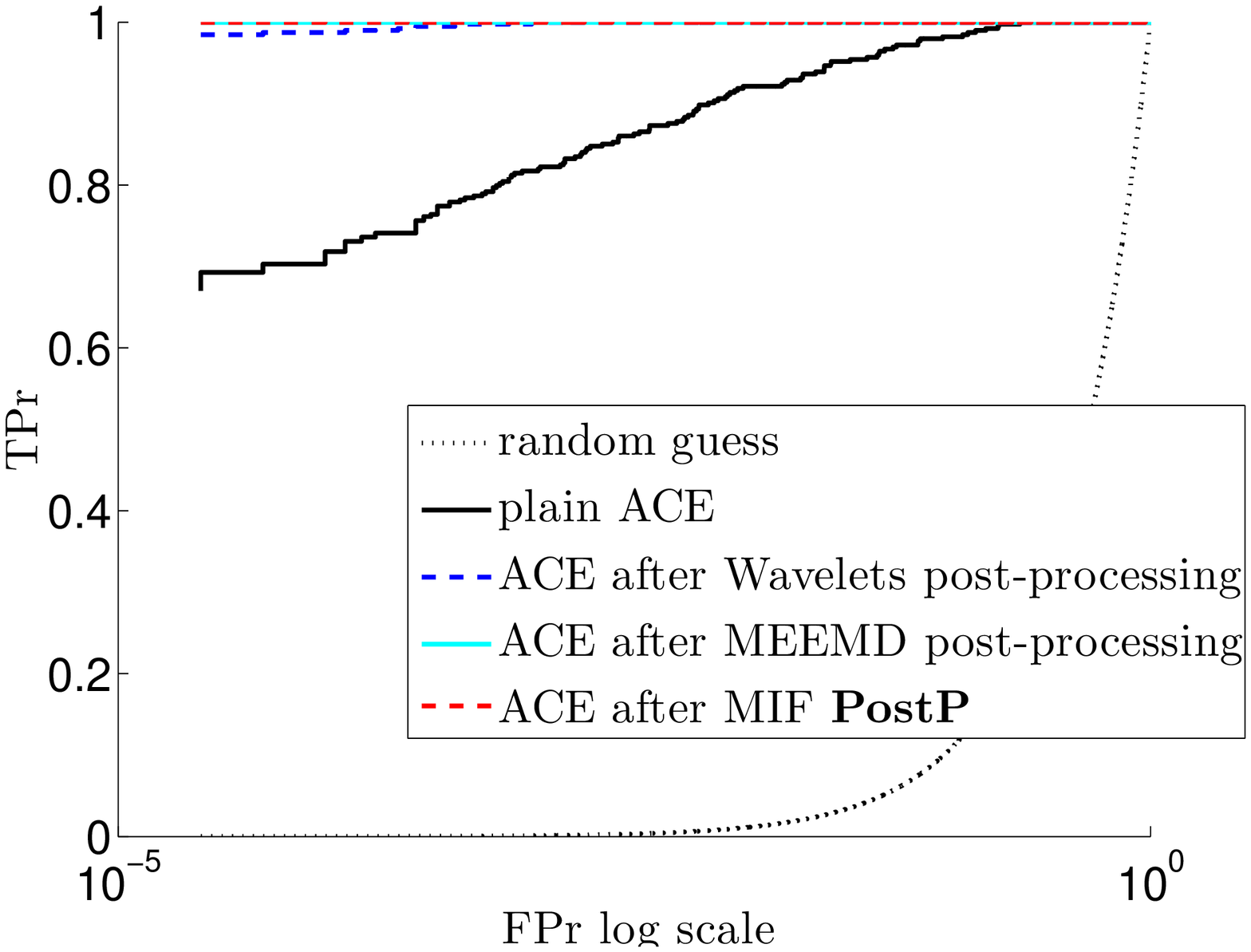}
     \caption{ROC curves  in log scale along the horizontal axis}
     \label{fig:dugway_rel_r134_roc_ACE_PP}
     \end{subfigure}
\caption{}\end{figure}

We start applying the ACE classifier to this hypercube. The outcomes are plotted in Figure \ref{fig:dugway_rel_r134_det_ACE} where we clearly see the contribution of the noise. The corresponding ROC curve is shown in Figure \ref{fig:dugway_rel_r134_roc_ACE_PP_noLog}, its AUC value is approximately $0.98919$. If we apply the post--processing method \textbf{PostP} described in the previous section we clean the classification preserving the shape of the plume as shown both in Figure \ref{fig:dugway_rel_r134_det_ACE_PP} and by the ROC curve plotted in Figure \ref{fig:dugway_rel_r134_roc_ACE_PP_noLog}. Its AUC value is approximately  $0.99508$. We compare this result with the  post--processing methods based on wavelet transform and on the Multi-dimensional empirical EMD (MEEMD) method, which is a multidimensional version of the ensemble EMD proposed by Huang et al. in \cite{wu2009multi}. For the first one we use the Matlab noise removal tool called \emph{Stationary Wavelet Transform Denoising 2--D}, where we choose the Daubechies db2 wavelet to get the cleaned classification shown in Figure \ref{fig:dugway_rel_r134_det_ACE_Wav_PP}, whose ROC curve is plotted in Figure \ref{fig:dugway_rel_r134_roc_ACE_PP_noLog}. The AUC value in this case is approximately $0.99507$. While for the MEEMD method we use the code provided in \cite{wu2009multi} where we set to 10 the number of perturbations. The post--processed classification is obtained subtracting the first 2 IMFs from the original classification. The corresponding cleaned classification is depicted in Figure \ref{fig:dugway_rel_r134_det_ACE_PP_MEEMD}, its ROC curve is shown in Figure \ref{fig:dugway_rel_r134_roc_ACE_PP_noLog} and the AUC value is exactly equal, up to machine precision, to the one obtained using the \textbf{PostP} approach.

 We point out that, since in all the examples under study the roc curves tend to be close each
other and the differences are concentrated in the top left corner of the ROC curve plots, from now
on we opt to plot ROC curves using a log scale along the horizontal axis, ref. Figure \ref{fig:dugway_rel_r134_roc_ACE_PP}.

Another observation regards the shape of the boundaries produced using these post--
processing techniques. As we already mentioned, the pixels close to the boundaries of the plume, shown in orange in Figure
\ref{fig:dugway_rel_r134_Plume}, are completely excluded in the computation of the ROC curves. This is the reason why all the
classifications end up having almost identical ROC curves and area under the curve values, even when they have pretty different pixel
classification values as for the case of the wavelet transform post--processing versus the one obtained with the \textbf{PostP} and the MEEMD techniques.

  To have a better understanding of the true performances of these different post--processing
methods we can compare directly the classifications depicted in Figures \ref{fig:dugway_rel_r134_det_ACE_PP},
\ref{fig:dugway_rel_r134_det_ACE_Wav_PP} and \ref{fig:dugway_rel_r134_det_ACE_PP_MEEMD} with the ground truth
shown in Figure \ref{fig:dugway_rel_r134_Plume}. From this comparison it is evident that when using the wavelet
transform, while the noise is removed in a reasonable way, the shape of the plume is partially
lost due to excessive blunting. Since in this context it is very important to identify the boundaries
of the plume with high accuracy, the wavelet transform proves to be not the best option as a
post--processing technique. More in general, in all the examples we tested, standard techniques
for noise removal, like Fourier or wavelets transform, tend to have problems producing accurate
boundaries. This is due to two reasons: this kind of signals are in general non-stationary, and
the standard transformations often use bases that are determined a priori. Both limitations are
overcome by the usage of the Multidimensional Iterative Filtering (MIF) method as a post--processing technique. In fact MIF technique can handle
a non-stationary signal and does not require any a priori assumption or knowledge on it, as
confirmed in this example by the shape of the plume and its boundary after MIF post--processing,
ref. Figure \ref{fig:dugway_rel_r134_det_ACE_PP}.

Regarding the results obtained via MEEMD we point out that, even if the results are comparable with the one produced with the proposed method, this approach has disadvantages. First of all it is still an open problem to prove the convergence of this technique even in one dimension. In fact in this case the method reduces to the so called EEMD technique, whose convergence is not known. Whereas we recall that, even if the convergence of the MIF method is still an open problem as well, its one dimensional version, the Iterative Filtering technique, has been proven to be convergent under mild conditions \cite{cicone2014adaptive}. Furthermore the MEEMD proves to be much slower than the proposed method. On an average personal computer it takes around 4 seconds to produce the clean classification plotted in Figure \ref{fig:dugway_rel_r134_det_ACE_PP}, while the MEEMD runs for more than 190 seconds on the same computer to produce the result shown in Figure \ref{fig:dugway_rel_r134_det_ACE_PP_MEEMD}.

We observe that with this dataset the MF classifier gives results similar to ACE, so we avoid presenting them and move directly to the COS classifier.

\begin{figure}[H]
\begin{subfigure}[b]{0.48\textwidth}
        \centering
        \includegraphics[width=\textwidth]{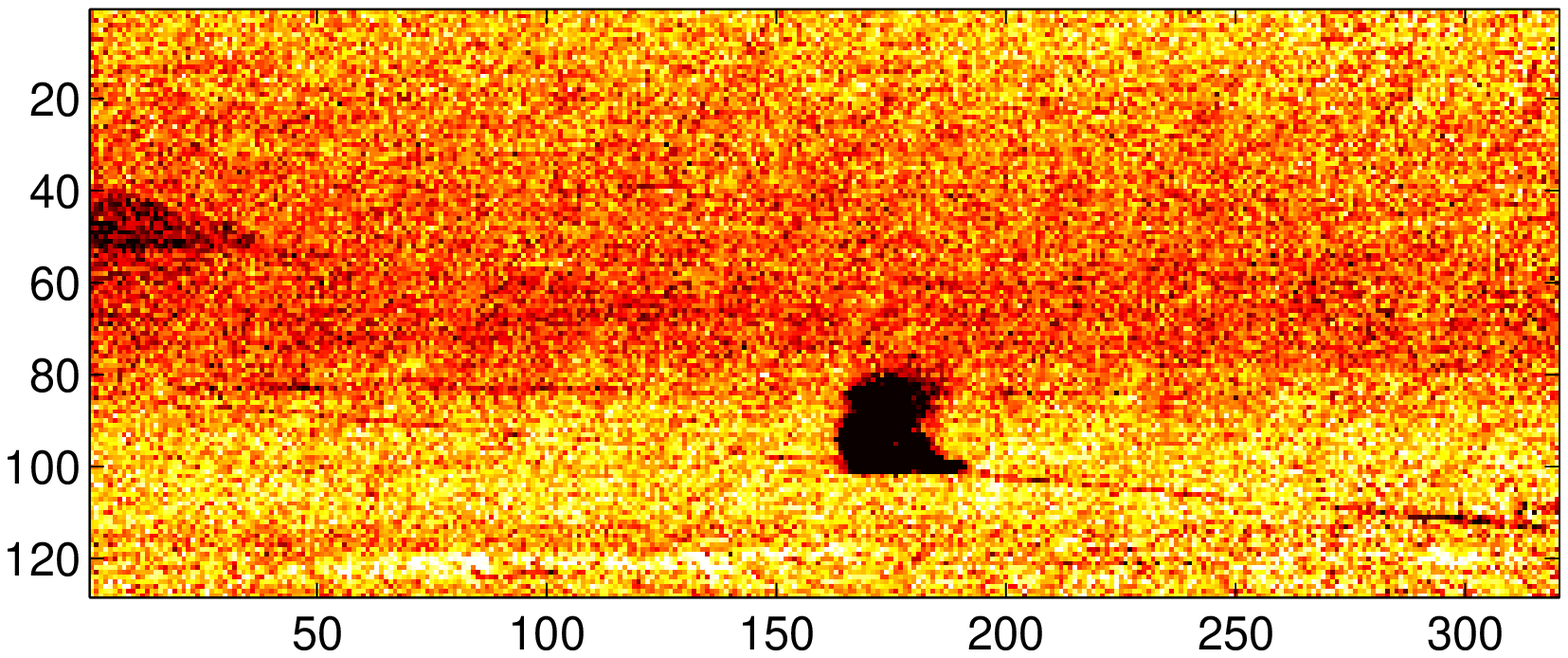}
        \caption{Detection map of the reversed COS classification of the raw data}
        \label{fig:dugway_rel_r134_det_1-COS}
        \end{subfigure}
      \hskip 2mm
     \begin{subfigure}[b]{0.48\textwidth}
     \centering
        \includegraphics[width=\textwidth]{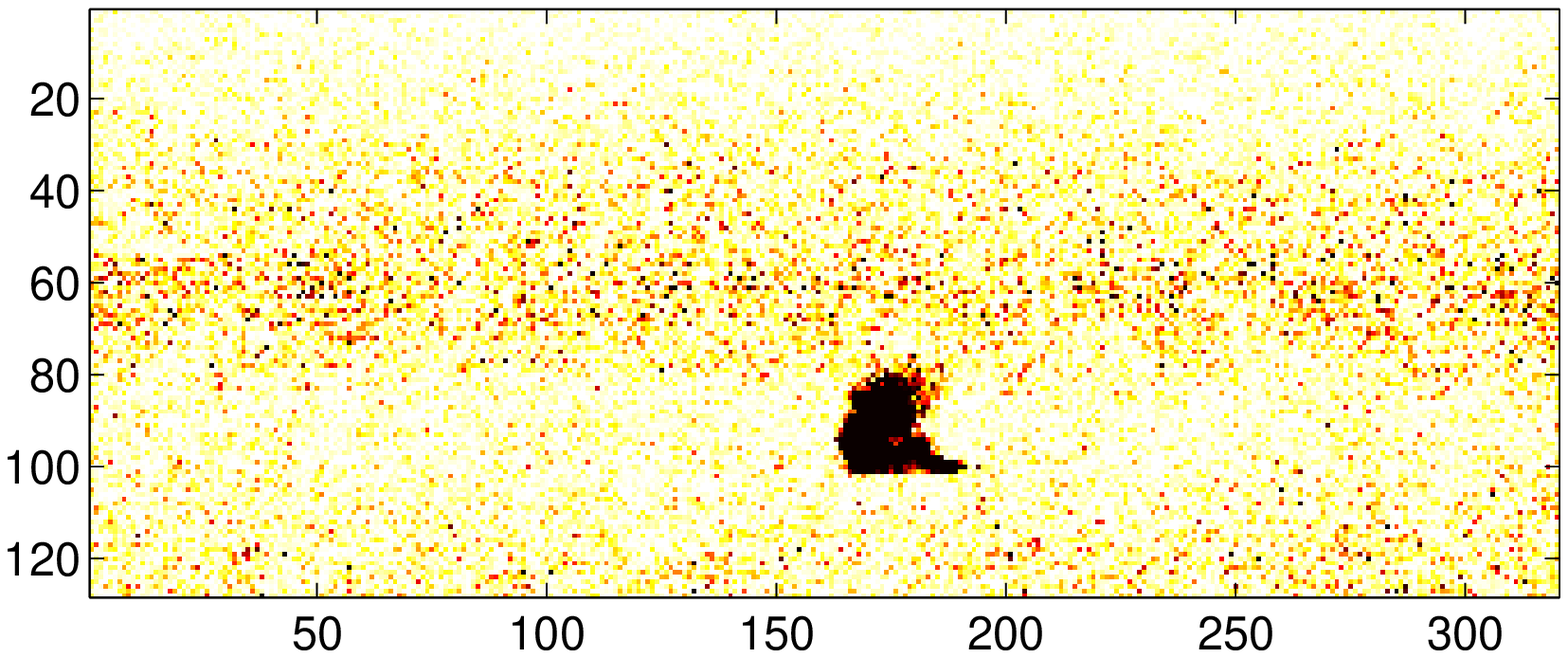}
                \caption{Detection map of the COS classification of the pre--processed data}
                \label{fig:dugway_rel_r134_det_pre_COS}
        \end{subfigure}
\caption{}\end{figure}
\begin{figure}[H]
\begin{subfigure}[b]{0.48\textwidth}
\centering
        \includegraphics[width=\textwidth]{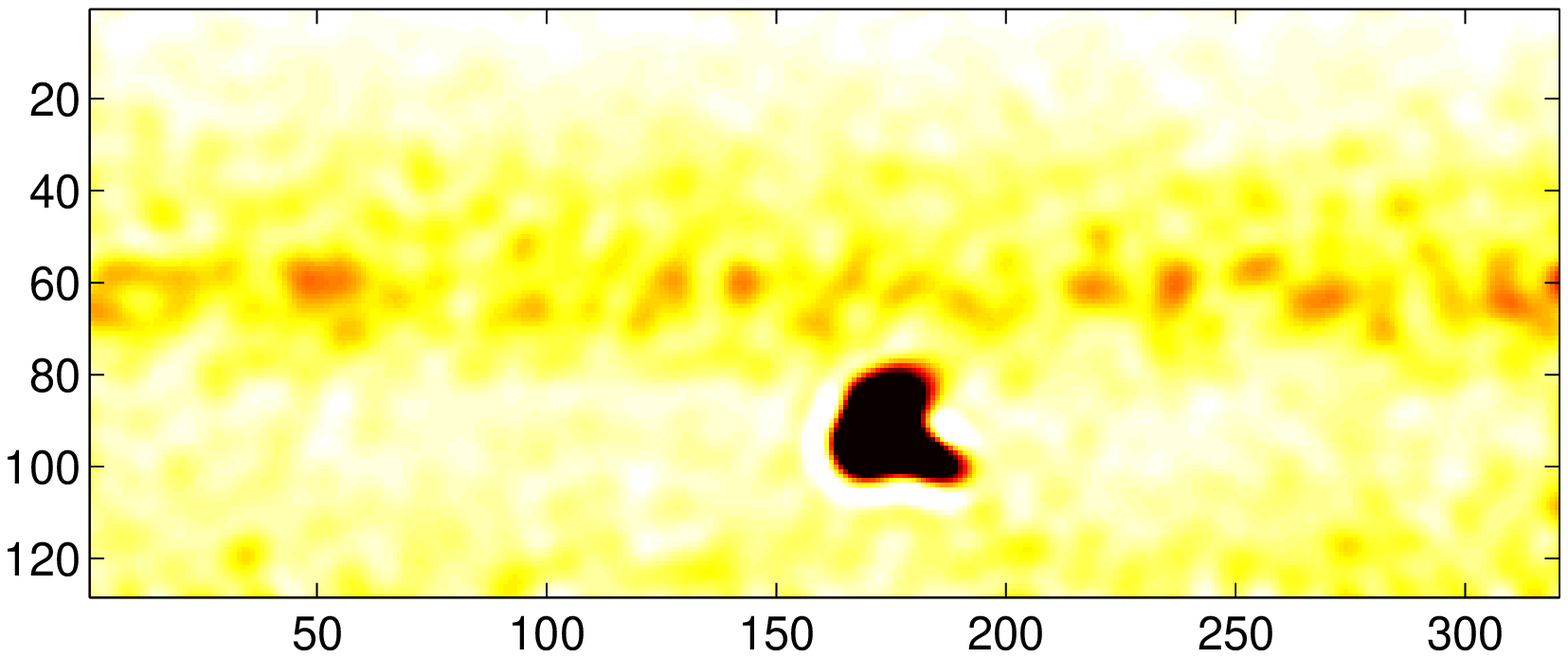}
                \caption{Detection map of the \textbf{PostP} post--processed COS classification of the pre--processed data}
                \label{fig:dugway_rel_r134_det_pre_COS_PP}
        \end{subfigure}
        \hskip 2mm
        \begin{subfigure}[b]{0.48\textwidth}
        \centering
        \includegraphics[width=\textwidth]{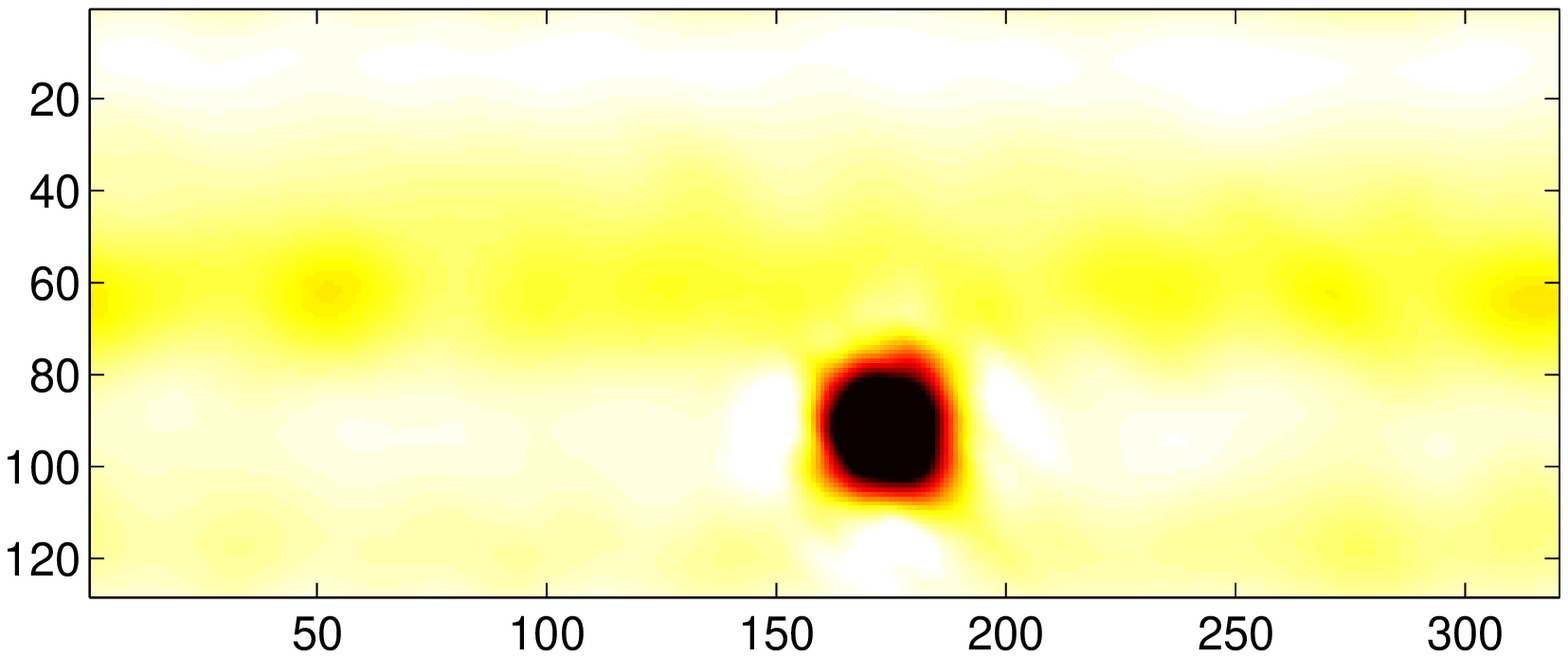}
        \caption{Detection map of the Wavelets post--processed COS classification of the pre--processed data}
        \label{fig:dugway_rel_r134_det_pre_COS_Wav_PP}
        \end{subfigure}
\caption{}
\end{figure}
\begin{figure}[H]
\begin{subfigure}[b]{0.48\textwidth}
     \centering
     \includegraphics[width=\textwidth]{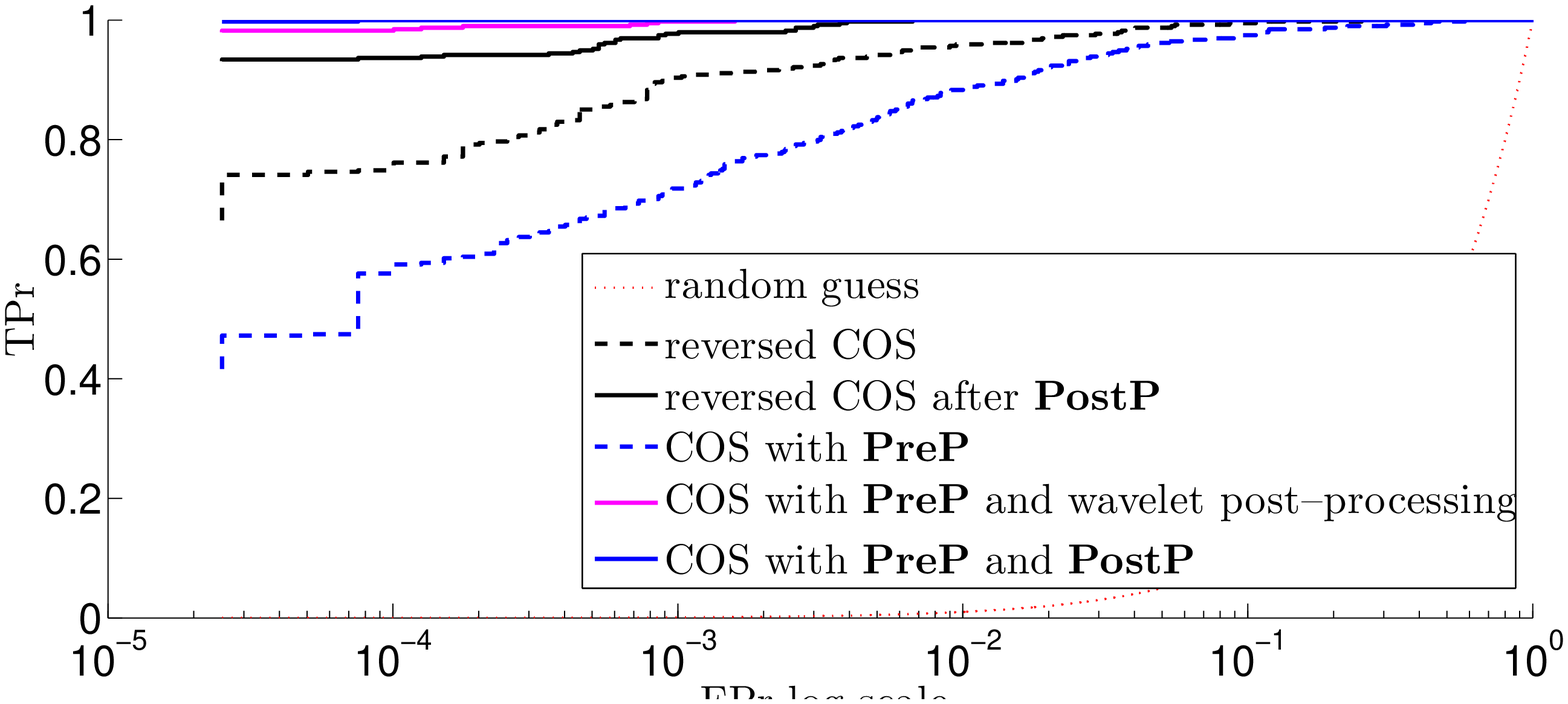}
     \caption{ROC curves  for the COS classification of the hypercube \emph{Location 1 released r134a}}
     \label{fig:dugway_rel_r134_roc_COS_PP}
     \end{subfigure}\hskip 2mm
\begin{subfigure}[b]{0.48\textwidth}
\centering
        \includegraphics[width=\textwidth]{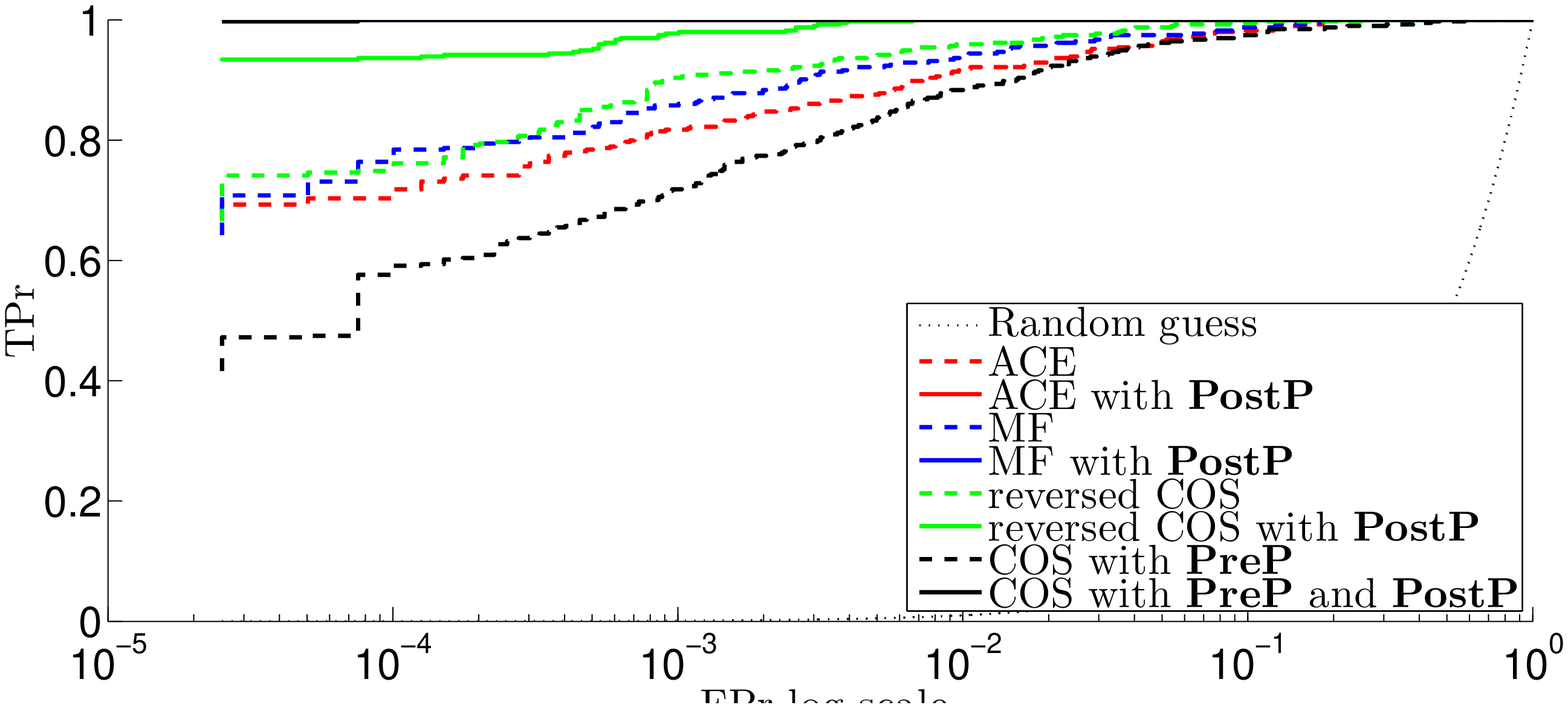}
                \caption{ROC curves  comparison for the hypercube \emph{Location 1 released r134a}}
                \label{fig:dugway_rel_r134_roc_curves}
        \end{subfigure}
\caption{}\end{figure}

The pixel classification produced by COS applied to the raw hypercube is shown in Figure \ref{fig:dugway_rel_r134_det_1-COS} where we reversed the pixel classification values. We need to do so in order to have meaningful results. In fact the corresponding ROC curve, plotted in dashed black in Figure \ref{fig:dugway_rel_r134_roc_COS_PP}, would be below the random guess curve without the reversion. If we apply the \textbf{PostP} method to the reversed COS classification we can improve the classification performance. This is testified by the increasing of the area under the ROC curve, ref. ROC curve plotted in solid black in Figure \ref{fig:dugway_rel_r134_roc_COS_PP}. However, even with post--processing, the COS classification of the raw hyperspectral image still needs to be reversed. To solve this issue we devise and apply the pre--processing method \textbf{PreP} based on MIF algorithm and described in the previous section. The COS classification of the pre--processed hypercube is shown in Figure \ref{fig:dugway_rel_r134_det_pre_COS}. Its ROC curve, plotted in dashed blue in Figure \ref{fig:dugway_rel_r134_roc_COS_PP}, is now directly above the random guess curve without any need of reversing the classification. Nevertheless its performances are worse than the ones produced by the reversed COS classification of the raw dataset, the dashed black curve is clearly above the solid blue one. However if we apply the \textbf{PostP} method to the COS classification of the pre--processed hypercube we get better performance than any other COS classification with or without post--processing. In fact the  solid blue ROC curve is above all the others in Figure \ref{fig:dugway_rel_r134_roc_COS_PP}. We can also try to apply the wavelet transform as a denoising post--processing method. If we plot the corresponding ROC curve, its performances are slightly worse than the ones obtained using the proposed method, the solid magenta curve is slightly lower than the solid blue one in Figure \ref{fig:dugway_rel_r134_roc_COS_PP}. Furthermore the boundaries end up being excessively blunted, Figure \ref{fig:dugway_rel_r134_det_pre_COS_Wav_PP}, for the reasons explained above, while the proposed post--processing technique allows to produce boundaries that look more natural, ref. Figure \ref{fig:dugway_rel_r134_det_pre_COS_PP}.

So, overall, the best performances are obtained from the combination of the pre-- and post--processing techniques, both based on MIF, applied to the COS classifier. We point out that the range of values produced by the COS classifier applied to the raw dataset is in the narrow interval $[0.8,\,0.84]$ after the reversion, while the one corresponding to the COS classification of the pre--processed dataset is in the wider interval $[0,\,0.5]$ without any need of reversing them. The performance of the COS classifier equipped with the pre--processing method becomes comparable with the one of ACE and MF algorithms.

A question that naturally arises at this point is why the pre--processing method works so well in enhancing the COS classification. The reason is evident from Figures \ref{fig:dugway_rel_r134_RawPlumePixel1} and \ref{fig:dugway_rel_r134_PreprocPlumePixel1}. Removing the global and local trend from the signature of each pixel allows to properly measure the angle between the pixels and the chemical signatures. Furthermore we observe that all the datasets provided for the aforementioned challenge always have the property that the signature of the chemical seems to be subtracted from the signature of the pixels contained in the plume. It remains an open problem to understand the reason of such a behavior of the signatures.
From this last observation it is evident that only after removing the trend of the pixel signatures the direct comparison of the signatures made by the COS classifier produces a meaningful classification. We observe that the combination of a pre-- and post--processing of the COS classification works well in increasing the performance of this classifier for all the datasets provided for the 2012 DTRA Chemical Detection Challenge.

\subsection{Blind case: hypercube ``Location 2 released sf6 blind'' }

\begin{figure}[H]
        \begin{subfigure}[b]{0.48\textwidth}
                \centering
                \includegraphics[width=\textwidth]{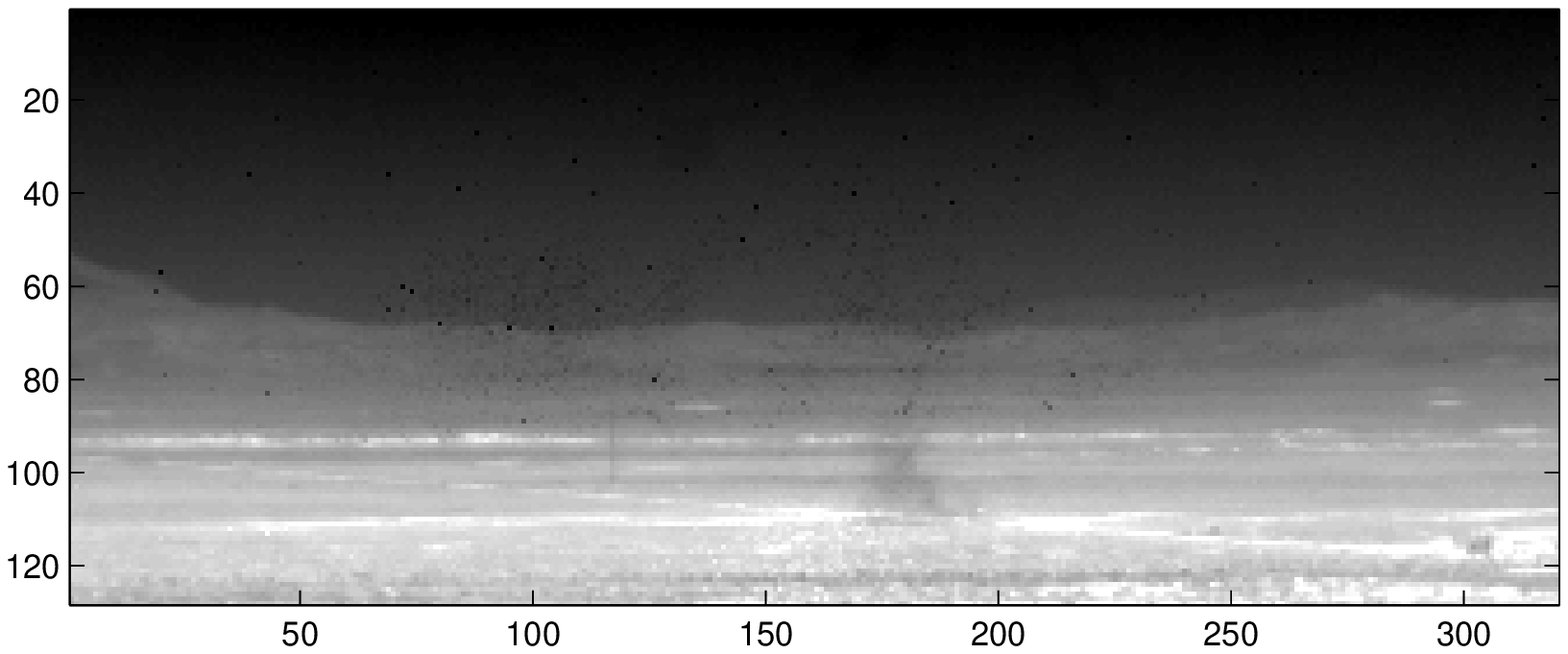}
                \caption{Contrast--enhanced spectral--mean image using the raw data \emph{Location 2 released sf6 blind}}
                \label{fig:dugway_rel_sf6_blind_Image}
        \end{subfigure}\hskip 2mm
        \begin{subfigure}[b]{0.48\textwidth}
        \centering
        \includegraphics[width=\textwidth]{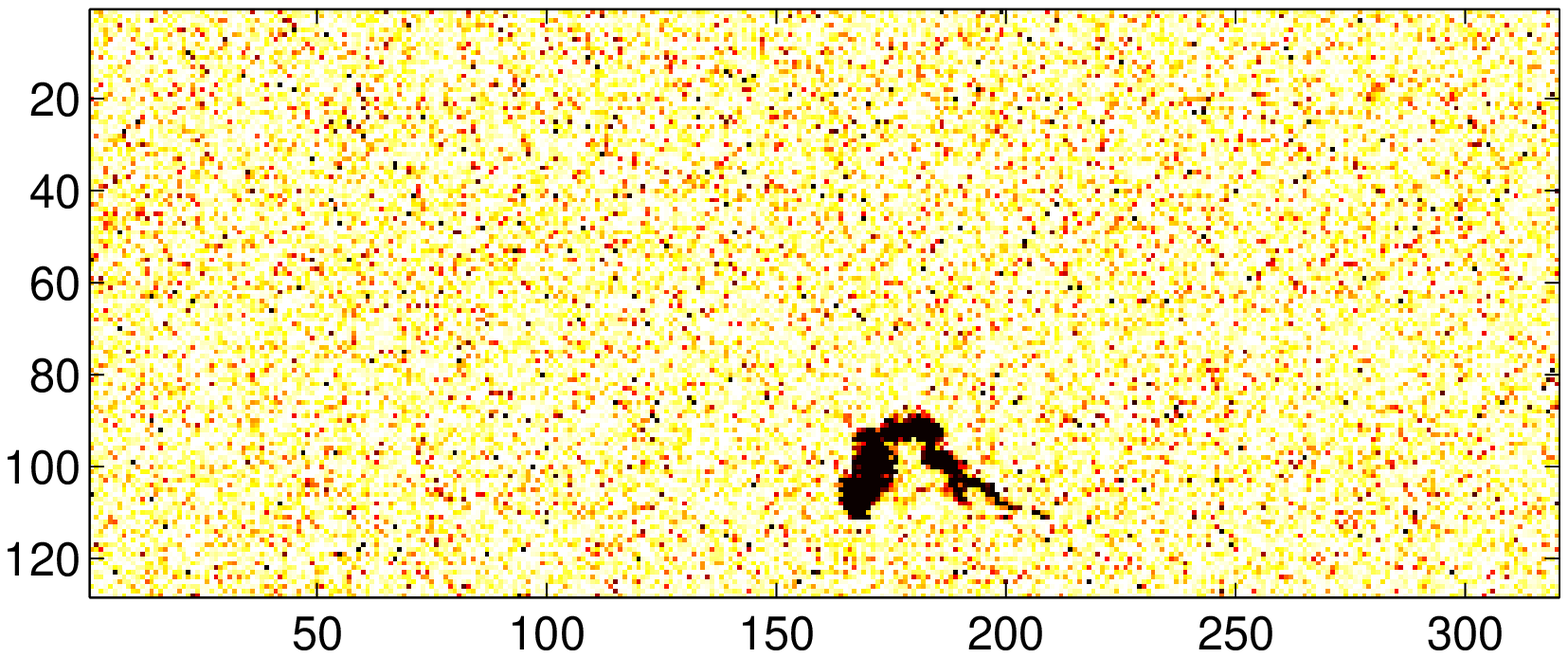}
        \caption{Detection map of the ACE classification of the raw data}
        \label{fig:dugway_rel_sf6_blind_det_ACE}
        \end{subfigure}
\caption{}\end{figure}

As another example we consider one of the blind datasets provided for the challenge. The hypercube image is shown in Figure \ref{fig:dugway_rel_sf6_blind_Image}. It is known that this image contains a plume of the sf6 chemical, whose signature is given.

\begin{figure}[H]
        \begin{subfigure}[b]{0.48\textwidth}
        \centering
        \includegraphics[width=\textwidth]{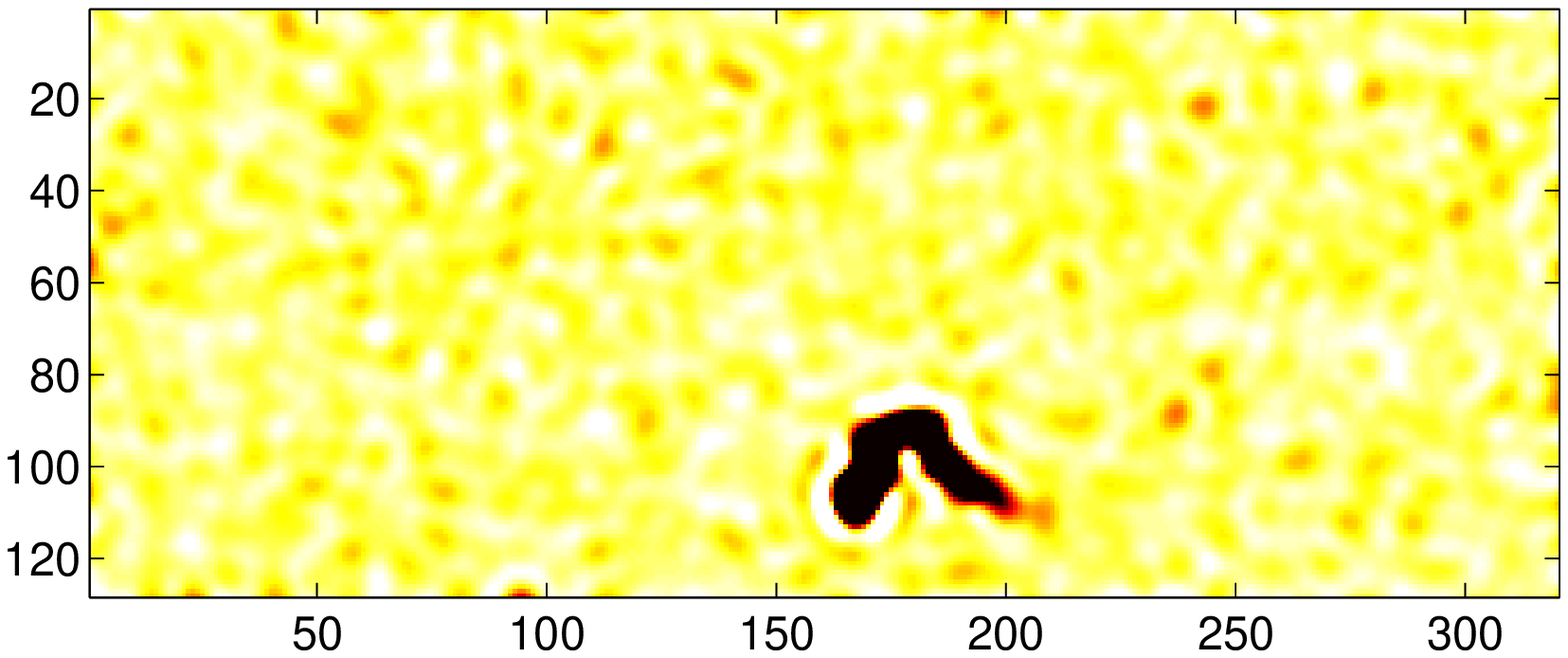}
                \caption{Detection map of the MIF post--processed ACE classification of the raw data}
                \label{fig:dugway_rel_sf6_blind_det_ACE_PP}
        \end{subfigure}\hskip 2mm
        \begin{subfigure}[b]{0.48\textwidth}
                \centering
                \includegraphics[width=\textwidth]{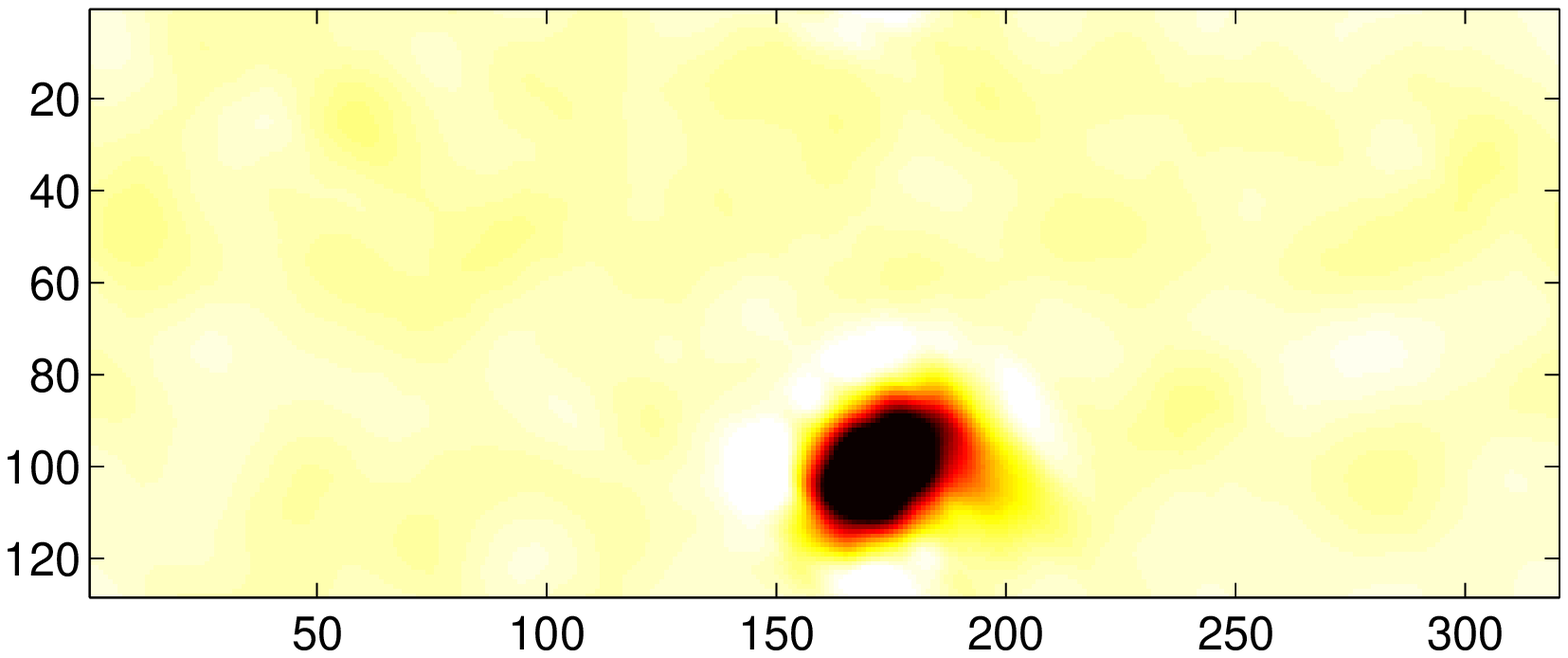}
        \caption{Detection map of the Wavelets post--processed ACE classification of the raw data}
        \label{fig:dugway_rel_sf6_blind_det_ACE_Wav_PP}
        \end{subfigure}
\caption{}\end{figure}

\begin{figure}[H]
        \begin{subfigure}[b]{0.48\textwidth}
        \centering
        \includegraphics[width=\textwidth]{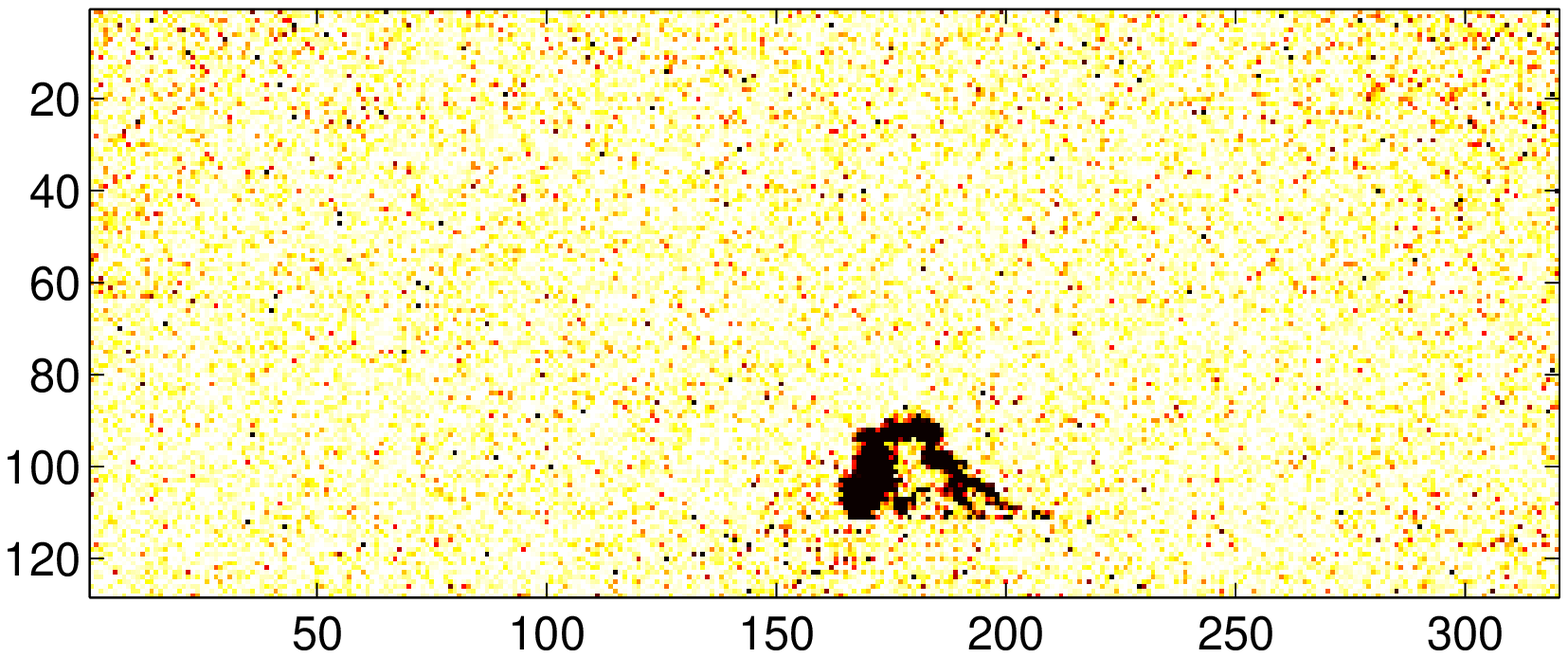}
        \caption{Detection map of the MF classification of the raw data}
        \label{fig:dugway_rel_sf6_blind_det_MF}
        \end{subfigure}\hskip 2mm
        \begin{subfigure}[b]{0.48\textwidth}
        \centering
        \includegraphics[width=\textwidth]{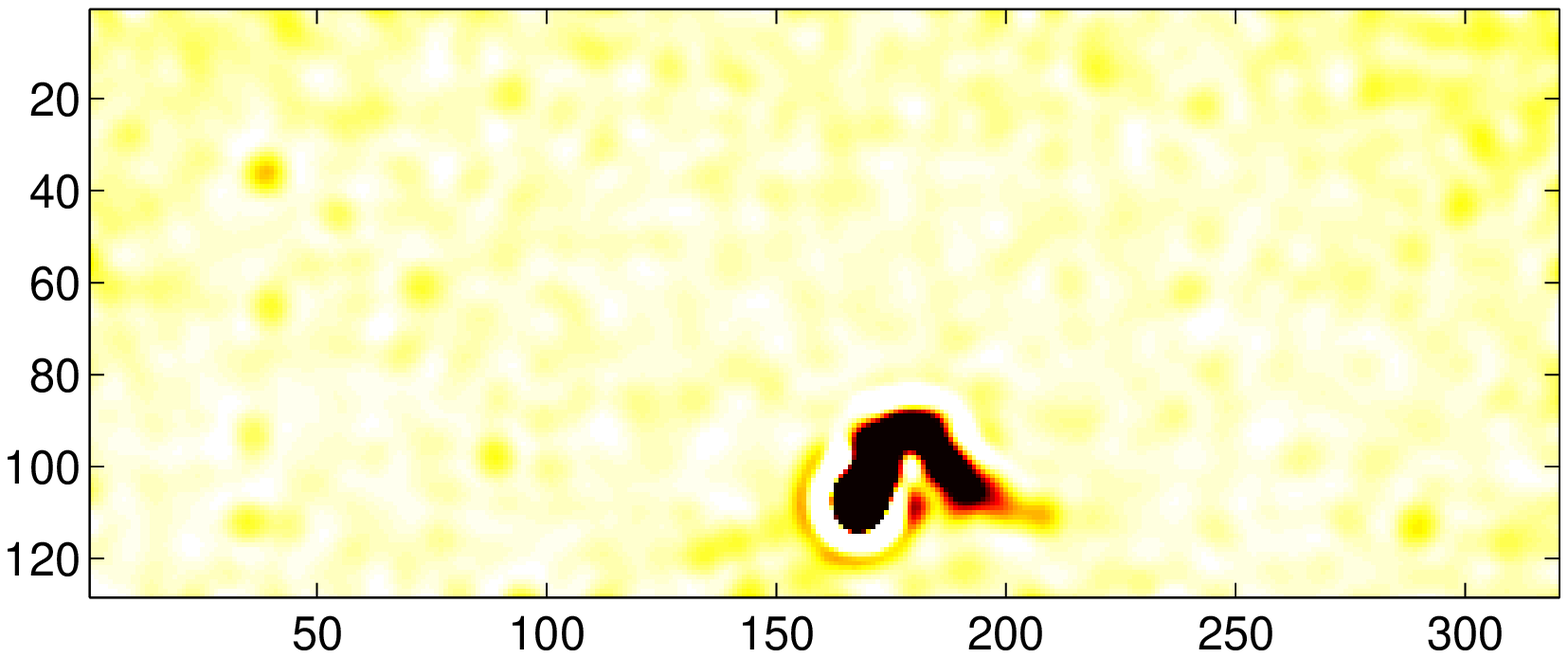}
                \caption{Detection map of the MIF post--processed MF classification of the raw data}
                \label{fig:dugway_rel_sf6_blind_det_MF_PP}
        \end{subfigure}
\caption{}\end{figure}

\begin{figure}[H]
        \begin{subfigure}[b]{0.48\textwidth}
                \centering
                \includegraphics[width=\textwidth]{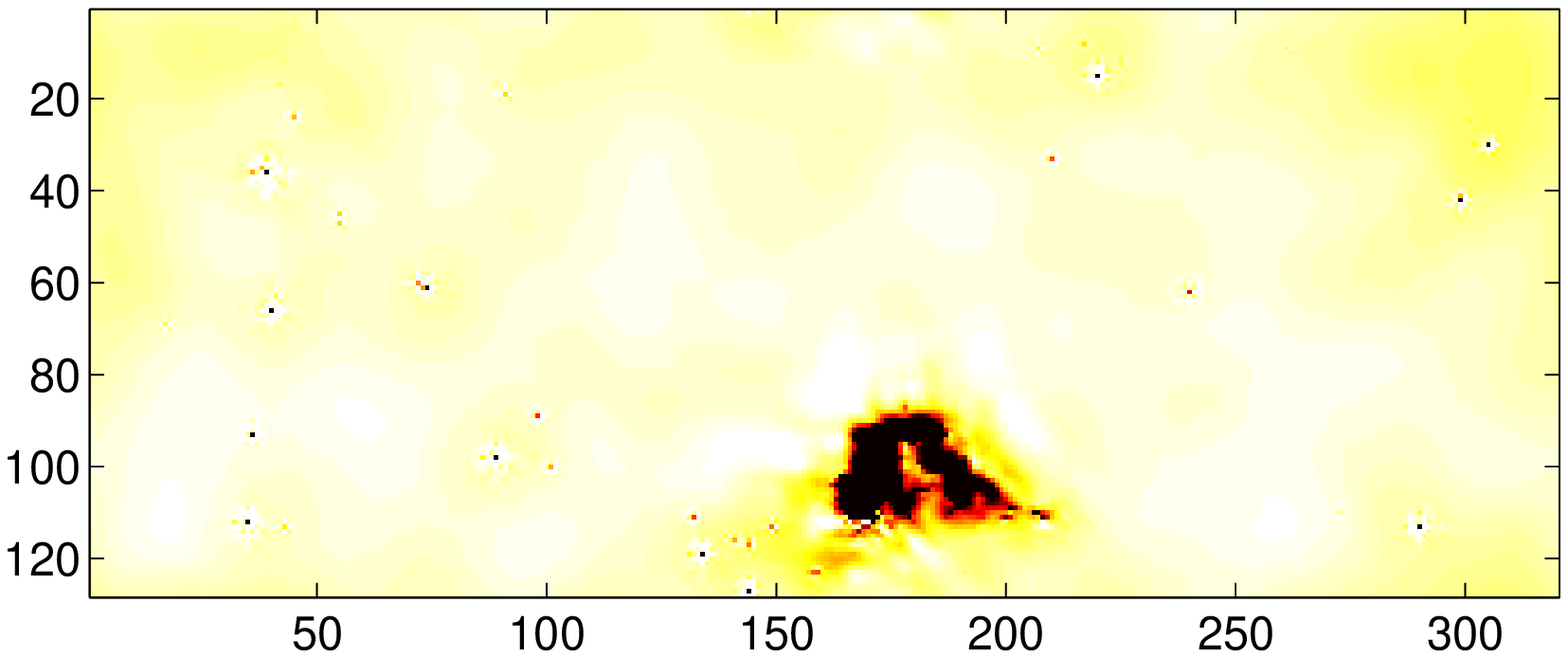}
        \caption{Detection map of the Wavelets post--processed MF classification of the raw data}
        \label{fig:dugway_rel_sf6_blind_det_MF_Wav_PP}
        \end{subfigure}\hskip 2mm
        \begin{subfigure}[b]{0.48\textwidth}
        \centering
        \includegraphics[width=\textwidth]{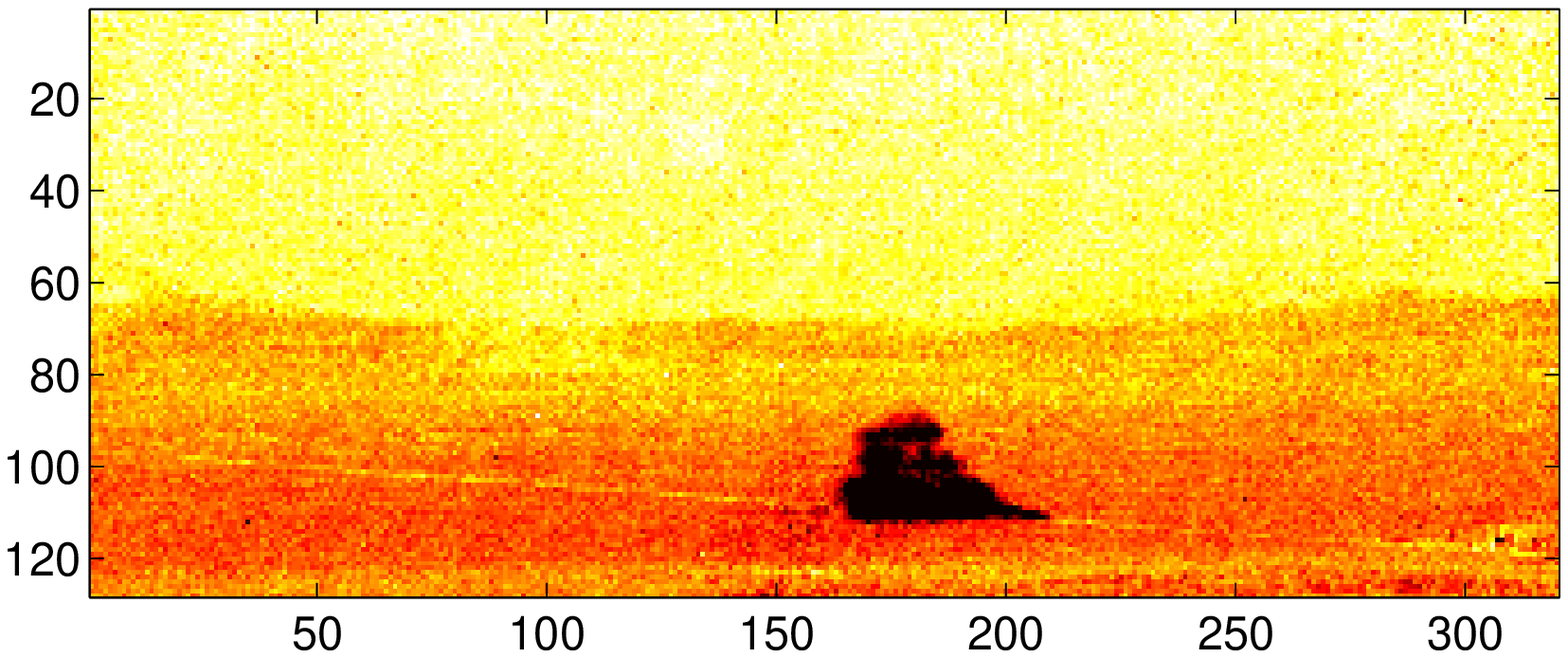}
        \caption{Detection map of the reversed COS classification of the raw data}
        \label{fig:dugway_rel_sf6_blind_det_1-COS}
        \end{subfigure}
\caption{}\end{figure}

If we apply the ACE classifier to the hypercube as it is, we obtain the classification depicted in Figure \ref{fig:dugway_rel_sf6_blind_det_ACE}. We can apply the \textbf{PostP} post--processing method to reenforce the spatial correlations, ref. Figure \ref{fig:dugway_rel_sf6_blind_det_ACE_PP}. The results produced by the Matlab noise removal tool called \emph{Stationary Wavelet Transform Denoising 2--D}, with the choice of the Daubechies db2 wavelet, are depicted in Figure \ref{fig:dugway_rel_sf6_blind_det_ACE_Wav_PP}. Similar results we get from the post--processing of the MF classification, ref. Figure \ref{fig:dugway_rel_sf6_blind_det_MF}. Using the proposed \textbf{PostP} technique, Figure \ref{fig:dugway_rel_sf6_blind_det_MF_PP}, we clean in a more effective way than with a standard technique, ref. Figure \ref{fig:dugway_rel_sf6_blind_det_MF_Wav_PP}.

\begin{figure}[H]
\begin{subfigure}[b]{0.48\textwidth}
\centering
        \includegraphics[width=\textwidth]{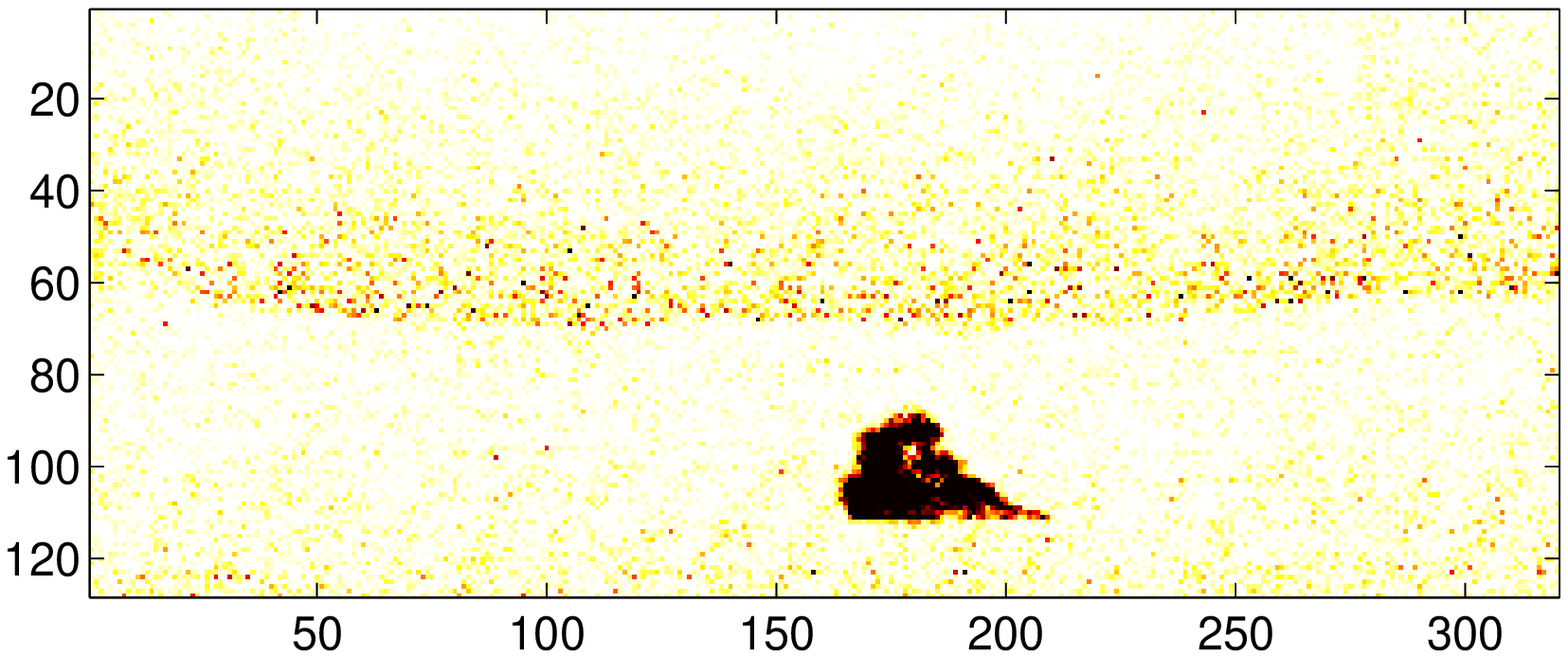}
                \caption{Detection map of the COS classification of the pre--processed data}
                \label{fig:dugway_rel_sf6_blind_det_pre_COS}
        \end{subfigure}\hskip 2mm
\begin{subfigure}[b]{0.48\textwidth}
\centering
        \includegraphics[width=\textwidth]{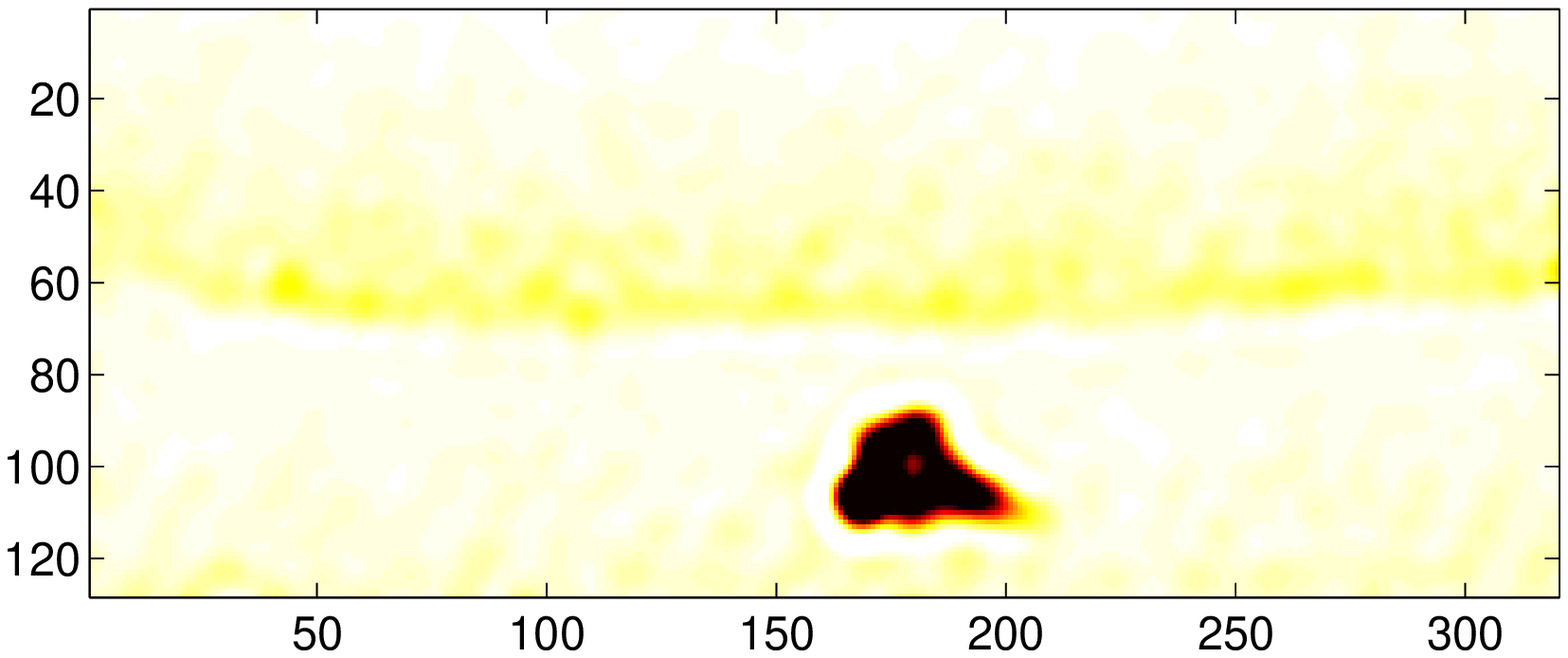}
                \caption{Detection map of the post--processed COS classification of the pre--processed data}
                \label{fig:dugway_rel_sf6_blind_det_pre_COS_PP}
        \end{subfigure}
\caption{}\end{figure}

With the plain COS classifier, as it was for the previous dataset,  we need to reverse the classification in order to produce meaningful results, Figure \ref{fig:dugway_rel_sf6_blind_det_1-COS}. If we apply the proposed  \textbf{PreP} pre--processing method to the hypercube we can remarkably increase the performance of the classifier as shown in Figure \ref{fig:dugway_rel_sf6_blind_det_pre_COS} so that there is no more need to reverse the classification and its performance become better than the one of the other classifier tested. In fact the COS method equipped with the proposed pre--processing method is producing a more natural plume shape than the one produced by any other classification techniques tested. So the pre--processing technique allows to transform a known failing classification method into what appears to be an extremely successful classification algorithm. We observe that the numerical range of the COS classification of to the raw hypercube spans approximately the interval $0.92$ and $0.94$, while the one of the pre--processed hypercube is in between 0 and $0.5$.

Also in the case of the COS classifier the application of the \textbf{PostP} post--processing technique helps in cleaning and improving the classification, Figure \ref{fig:dugway_rel_sf6_blind_det_pre_COS_PP}.

\section{Conclusions}

Inspired by the DTRA 2012 Chemical--Detection Algorithm Challenge posted by Dimitris Manolakis and his research group at MIT Lincoln Laboratory, in this paper we tackle the problem of the accurate detection of the boundaries of a chemical plume hidden in a hyperspectral image by means of  classification algorithms. In particular we devise a post-- and a pre--processing algorithm aimed to improve the performance of known classifiers.

The post--processing technique we present in this paper is based on the Multidimensional Iterative Filtering method \cite{cicone2015multi} which is an algorithm for the decomposition of non--stationary signals. This post--processing technique is a uniform procedure which allows to improve the performance of any classifier used in the chemical plume detection problem by means of an adaptive and data driven cleaning and smoothing of its classification. The choice of the MIF method appears to be ideal, since it is an algorithm developed to handle naturally non--stationary signals, which is the case for almost any real life signal. Furthermore MIF algorithm does not require to make any a priori assumption on the kind of data we are analyzing as for other techniques like wavelet transform and similar methods. We test this approach on the datasets provided for the aforementioned challenge, and we show some of these results in the Examples section.

Secondly we propose a pre--processing technique which, independently from the dataset under study, allows to
decorrelate and mean--center the data. The goal of this procedure is to make the score values produced by the COS classifier equal to the ones produced by the ACE classifier, as shown in \eqref{eq:ACEeqCOS}. This is made possible thanks to the employment of the Multidimensional Iterative Filtering method which allows to decompose a signal, in particular a non--stationary one, in an adaptive and data driven way. From the tests run, and partially shown in the examples section, we notice that the chemical plume boundaries detected using the COS classifier equipped with the \textbf{PreP} pre--processing method appear to be meaningful. These results are confirmed by the ROC curves and the shape of the detected plumes. In some cases the detected boundaries prove to be more meaningful and reasonable than the ones produced using other standard classification methods.
In conclusion, given the Cosine Similarity measure, which is known to be a failing classification technique, the proposed pre--processing technique allows to transform it into a successful classifier which can even outperform in some cases other commonly used classifiers.

A possible direction for future researches is in the analysis of plumes containing mixtures of two or more known chemicals. Traditionally such mixtures have been assumed to be linear, however it is known that such assumption is far from reality. New ways for tackling this problem have to be devised. We point out also that, from an algorithmic point of view, the convergence of the Multidimensional Iterative Filtering method has been proved only for the 1D case \cite{cicone2014adaptive}. For higher dimensions its convergence has to be studied yet.

\section{Acknowledgments}

This work was supported by NSF Faculty Early Career Development (CAREER) Award DMS--0645266, DMS--1042998, DMS--1419027, and ONR Award N000141310408.

Antonio Cicone acknowledges support by National Group for Scientific Computation (GNCS -- INdAM) ``Progetto giovani ricercatori 2014'', and by Istituto Nazionale di Alta Matematica (INdAM) ``INdAM Fellowships in Mathematics and/or Applications cofunded by Marie Curie Actions''.


\end{document}